\documentclass[11pt]{article}
\usepackage{amsmath,amsfonts}
\usepackage{xcolor}
\usepackage{amsthm}         
\usepackage{dsfont}			
\usepackage{hyperref}
\usepackage{multirow}
\usepackage{graphicx}
\usepackage{subcaption}
\usepackage{epsfig}
\def\cb{{\bf c}}
\def\tb{{\bf t}}
\def\Xb{{\boldsymbol X}}
\def\Yb{{\boldsymbol Y}}

\newcommand*{\Ztu}{{}_{s}\mathring{\mathbb{Z}}_u}
\newcommand*{\thetab}{{\boldsymbol{\theta}}}
\newcommand*{\etab}{{\boldsymbol{\eta}}}

\newcommand*{\lb}{{\boldsymbol{\ell}}}
\newcommand*{\E}[1]{\operatorname{E}_{#1}}
\newcommand*{\var}[1]{\operatorname{var}_{#1}}
\newcommand*{\Prob}[2]{\operatorname{P}_{#1} \left( {#2} \right)}
\newcommand*{\Po}[1]{\operatorname{P}_{#1}}
\newcommand*{\deq}{{\, \stackrel{\text{def}}{=} \; }}
\newtheorem{theorem}{Theorem}
\newtheorem{remark}{Remark}
\newtheorem{lemma}[theorem]{Lemma}

\hypersetup{
	pdfborder = {0 0 0},
	urlbordercolor = {0 0 0},
	colorlinks = {false},
}

\newcommand\smallO{
	{{\scriptscriptstyle\mathcal{O}}}
}
\def\underacc #1/#2{\mathchoice{\uacc\textstyle{#1}{#2}}{\uacc\textstyle{#1}{#2}}
	{\uacc\scriptstyle{#1}{#2}}{\uacc\scriptscriptstyle{#1}{#2}}}
\def\uacc#1#2#3{\mathop{#2{}}\limits_{#1#3{}}}

\newcommand*{\ind}[1]{\mathds{1} \left( {#1} \right)}

\begin{document}
	
	\title{Exact variable selection in sparse nonparametric models}
	\author{Natalia Stepanova$^a$
		and Marie Turcicova$^{b0}$\\
		\\
		\small {\textit{$^a$School of Mathematics and Statistics, Carleton University, 4302 Herzberg Laboratories,}}\\
		\small {\textit{1125 Colonel By Drive, Ottawa, ON, K1S 5B6, Canada} } \\ 
		\small{\textit{$^b$Institute of Atmospheric Physics, Czech Academy of Sciences, }} \\
		\small{\textit{Bocni II 1401, Prague, 141 31,  Czech Republic}}}
	\date{} 
	
	\footnotetext{Corresponding author. Email address: mt@ufa.cas.cz}
	\maketitle

\begin{abstract}
We study the problem of adaptive variable selection in a Gaussian white noise model of intensity $\varepsilon$ under certain sparsity and regularity conditions on
an unknown regression function~$f$.
The $d$-variate regression function $f$ is assumed to be a sum of functions each depending on a smaller number $k$ of variables  ($1\leq k\leq d$). These functions are unknown to us and only few of them are nonzero. We assume that $d=d_\varepsilon\to \infty$ as $\varepsilon\to 0$ and consider the cases when $k$  is fixed and when $k=k_\varepsilon\to \infty$, $k=\smallO(d)$ as $\varepsilon\to 0$.
In this work, we introduce an adaptive selection procedure
that, under some model assumptions, identifies exactly all nonzero $k$-variate components of $f$.
In addition, we establish conditions under which exact identification of the nonzero components is impossible. These conditions ensure that
the proposed selection procedure is the best possible in the asymptotically minimax sense with respect to the Hamming risk.
\\ 
%\\ \smallskip

\textit{2020MSC}: Primary: 62G08, Secondary: 62H12, 62G20
\end{abstract}

%%%%%%%%%%%%%%%%%%%%%%%%%%%%%%%%%%%%%%%%%%%%%%
%% Please use \tableofcontents for articles %%
%% with 50 pages and more                   %%
%%%%%%%%%%%%%%%%%%%%%%%%%%%%%%%%%%%%%%%%%%%%%%
%\tableofcontents

\section{Introduction}\label{sec:problem_statement}

Mathematically, many common ``regular'' regression models with normally distributed errors can be accommodated within one functional model,
known as the Gaussian  white noise model.
Suppose that an unknown function $f(\tb)$, $\tb\in[0,1]^d$, of $d$ variables
is observed in the Gaussian white noise model
\begin{gather}\label{model1}
	X_{\varepsilon}=f+\varepsilon W,
\end{gather}
where $W$ is a $d$-dimensional Gaussian white noise, $\varepsilon>0$ is the noise intensity, and
$f\in L_2([0,1]^d)=L_2^d$ is a regular enough function.
We are interested in the case when $d=d_\varepsilon\to \infty$ as $\varepsilon\to 0.$
{In this model, an ``observation'' $X_{\varepsilon}$ is a~function $X_\varepsilon: L_2^d\to {\cal G}$ taking its values on the set ${\cal G}$
	of normal random variables such that if $\xi=X_{\varepsilon}(\phi)$ and $\eta=X_{\varepsilon}(\psi)$, where $\phi,\psi\in L_2^d$,
	then ${\rm E}(\xi)=(f,\phi)_{L_2^d}$, ${\rm E}(\eta)=(f,\psi)_{L_2^d}$, and ${\rm cov}(\xi,\eta)=\varepsilon^2(\phi,\psi)_{L_2^d}$,
	with $( \cdot,\cdot)_{L_2^d}$ being the inner product in $L_2^d$. For any $f\in L_2^d$, the observation $X_\varepsilon$
	determines the Gaussian measure ${\rm P}_{\varepsilon,f}$ on the Hilbert space $L_2^d$ with mean function $f$ and covariance operator
	$\varepsilon^2I$, where $I$ is the identity operator (for references, see \cite{GN2016}, \cite{IH.97}, and \cite{SK}).}
We assume, in addition to regularity constraints, that $f$ has a sparse structure, and
consider the problem of recovering the sparsity pattern of $f$
by using the asymptotically minimax approach.

In high-dimensional inference problems, including the problem under study, the curse of dimensionality phenomenon is unavoidable and is
reflected in the optimal rates of convergence
which are generally much slower as compared to those for one- or low-dimensional problems.
To deal with the curse of dimensionality in high-dimensional nonparametric problems, the functional ANOVA models,
which effectively reduce the ``working dimension'' of the problem, are often used.
One common functional ANOVA model assumes that
the $d$-variate function $f(\tb)$  observed in the Gaussian white noise model can be expressed as (see, for example, \cite{IS-2015} and \cite{Lin-2000})
\begin{gather*}
	f(\tb)=\sum_{u\subseteq\{1,\ldots,d\}} f_u(\tb_u),\quad \tb\in[0,1]^d,
\end{gather*}
where $\tb_{u}$ is the $\#(u)$-dimensional vector ($\# (\cdot)$ stands for cardinality) containing those components
of~$\tb$, whose indices belong to $u$, $f_u={\rm const}$ if $u=\emptyset$,  and
\begin{gather}\label{orthcon}
	\int_{0}^1 f_u(\tb_u)\,dt_j=0,\quad \mbox{for}\; j\in u,
\end{gather}
that is, the terms $f_u$ are mutually
orthogonal. Each function $f_u(\tb_u)$ depends only on variables in $\tb_u$ and describes the ``interaction''
between these variables.

Formalizing the notion of sparsity adopted in this article, we assume that the regression function~$f$ in model (\ref{model1})
admits the following sparse ANOVA-type decomposition.
For $1\leq k\leq d$, let ${\cal U}_{k,d}$ be the set of all subsets $u\subseteq\{1,\ldots,d\}$ of cardinality $k$, that is,
\begin{gather*} \label{def:Ukd}
	{\cal U}_{k,d}=\{u:u\subseteq\{1,\ldots,d\}, \#(u)=k\}.
\end{gather*}
Note that $\#\left( \mathcal U_{k,d}\right)={d\choose k}.$
If $u=\{j_1,\ldots,j_k\}\in {\cal U}_{k,d}$, $1\leq j_1<\ldots<j_k\leq d$, we denote as before $\tb_{u}=(t_{j_1},\ldots,t_{j_k})\in[0,1]^k$ and
assume that
\begin{gather}\label{f}
	f(\tb)=\sum_{u\in{\cal U}_{k,d}}\eta_u f_u(\tb_u),\quad \tb\in[0,1]^d,
\end{gather}
where $f_u$ satisfies (\ref{orthcon}) for $u\in{\cal U}_{k,d}$,
and the $\eta_j$s are unknown but deterministic quantities taking values in $\{0,1\}$:
$\eta_u=0$ means that the component $f_u$ is inactive, whereas $\eta_u=1$ means that the component $f_u$ is active.
The number $\sum_{u\in{\cal U}_{k,d}  }\eta_u$ of active components
is set to be small compared to the total number of components ${d\choose k}$, specifically
(recall that $d=d_\varepsilon\to \infty$ as $\varepsilon\to 0$)
\begin{gather}\label{sparsitycond}
	\sum_{u\in{\cal U}_{k,d}}\eta_u={d\choose k}^{1-\beta}(1+\smallO(1)),\quad \varepsilon\to 0,
\end{gather}
where $\beta\in(0,1)$ is the \textit{sparsity index}. We may think of $\sum_{u\in{\cal U}_{k,d}  }\eta_u$
as the integer part of ${d\choose k}^{1-\beta}$, denoted by $ \left\lfloor{d\choose k}^{1-\beta}\right\rfloor$, and introduce the set
\begin{equation*}
	\mathcal{H}^d_{k,\beta} = \left\{ \boldsymbol{\eta}=( \eta_u)_{ u \in \mathcal{U}_{k,d}}: \eta_{u} \in \{0,1\}\;\mbox{and condition}\;
	(\ref{sparsitycond})\;\mbox{holds}\right\}.
	%\label{def:H_dbk}
\end{equation*}
The orthogonality conditions in (\ref{orthcon}) imply that if  $u\neq v$ are subsets of ${\cal U}_{k,d}$, then $f_u(\tb_u)$ and $f_{v}(\tb_v)$
are orthogonal (in $L_2^d$) to each other and to a constant, which guarantees uniqueness of representation~(\ref{f}).
The signal $f$ as in (\ref{f}) is \textit{sparse} because the majority of the components $f_u$  are inactive,
and only $\left\lfloor{d\choose k}^{1-\beta}\right\rfloor$ components are active,
where ${d\choose k}^{1-\beta}=\smallO\left( {d\choose k}\right)$ as $d\to \infty$ and $k=\smallO(d)$.
In other words,  $f$~is the sum of a small number of $k$-variate functions.
Sparse models of this type with $k=1$, that is, sparse additive Gaussian white noise models,
have been studied in a number of publications dealing with high-dimensional inference problems
(see, for example,  \cite{BSt-2017}, \cite{DSA-2012}, \cite{IS-2015}).

Based on the ``observation'' $X_\varepsilon$ in model (\ref{model1}) to (\ref{sparsitycond}), we wish to identify,
with high degree of accuracy, the active (nonzero) components $f_u$ of $f$, that is,
we wish to construct a good estimator $\hat{\etab}=\hat{\etab}(X_\varepsilon)=(\hat{\eta}_u)_{u\in {\cal U}_{k,d}}$ of $\etab=(\eta_u)_{u\in {\cal U}_{k,d}}$
taking values in $\{0,1\}^{{d\choose k}}$
that would tell us which terms $f_u$ in the sparse decomposition (\ref{f}) are active. This may be viewed as a variable selection
problem, and $\hat{\etab}$ may be named a selector.

In order to better situate this work in current literature, consider a simpler
vector model
\begin{gather*}
	X_u=\mu_d\eta_u+\varepsilon_u,\quad u\in {\cal U}_{k,d},
\end{gather*}
where the errors $\varepsilon_u$ are iid standard normal random variables, the quantities $\eta_u$ are as before and, in particular, condition (\ref{sparsitycond})
holds,
and $\mu_d>0$ is  the signal. When $\eta_u=0$, the observation $X_u$ is pure noise, whereas
when  $\eta_u=1$, it is the signal $\mu_d$ plus noise. Consider the recovery of $\etab=(\eta_u)_{u\in {\cal U}_{k,d}}$ in this model.
To quantify the performance  of a selection procedure that identifies nonzero components of $\etab$, the Hamming risk $\operatorname{E}_{\mu_d,\etab} |\hat{\etab} - \etab| = \operatorname{E}_{\mu_d,\etab} \left( \sum_u |\hat{\eta}_u - \eta_u| \right)$ can be used. In the asymptotic setup as $d\to \infty$, the recovery of $\etab=(\eta_u)_{u\in {\cal U}_{k,d}}$ is meaningful only when $\mu_d\to \infty$ as $d\to \infty$.
Indeed, define a~selector ${\etab}^*=({\eta}^*_u)_{u\in {\cal U}_{k,d}}$ as follows:
\begin{equation}
	{\eta}^*_u={\eta}^*_u(X_u) = \ind{ X_{u} > \sqrt{(2+\kappa_d) \log {d\choose k} }},\quad u\in {\cal U}_{k,d},
	\label{vector_selector}
\end{equation}
where $\kappa_d>0$ is such that $\kappa_d \to 0$ and $\kappa_d\log {d\choose k}\to \infty$ as $d\to \infty$.
One can show that  if
\begin{equation}
	\liminf_{d \to \infty} \frac{\mu_d}{\sqrt{\log {d\choose k}}} > \sqrt{2} (1+\sqrt{1-\beta}),
	\label{vector_cond:inf}
\end{equation}
then the selector ${\etab}^*=({\eta}^*_u)_{u\in {\cal U}_{k,d}}$  provides \textit{exact recovery} of
the sparsity pattern in the sense that
\begin{gather*}
	\limsup_{d\to \infty}\sup_{\etab\in H_{d,\beta}} \operatorname{E}_{\mu_d,\etab}\left(\sum_{u\in {\cal U}_{k,d}}|{\eta}^*_u-\eta_u|\right)=0,
\end{gather*}
where $H_{d,\beta}=\Big\{\etab=(\eta_u)_{u\in {\cal U}_{k,d}}: \eta_u\in \{0,1\}\;\mbox{and}\; \sum_{u\in {\cal U}_{k,d}}\eta_u\leq c{d\choose k}^{1-\beta} \;\mbox{for some } c>0\Big\}.$
In words, under condition (\ref{vector_cond:inf}),  the maximum Hamming risk of ${\etab}^*$ tends to zero as $d$ tends to infinity, and hence
${\etab}^*$ identifies correctly all the positions $u$ where the signal appears. At the same time, if
\begin{equation}
	\limsup_{d \to \infty} \frac{\mu_d}{\sqrt{\log {d\choose k}}} < \sqrt{2} (1+\sqrt{1-\beta}),
	\label{vector_cond:inf2}
\end{equation}
then
\begin{gather*}
	\liminf_{d\to \infty}\inf_{\tilde{\etab}}\sup_{\etab\in H_{d,\beta}} \operatorname{E}_{\mu_d,\etab}\left(\sum_{u\in {\cal U}_{k,d}}|\tilde{\eta}_u-\eta_u|\right)>0,
\end{gather*}
where the infimum is over all selectors $\tilde{\etab}=(\tilde{\eta}_u)_{u\in{\cal U}_{k,d}}$, and hence exact variable selection is impossible.
The above stated results hold true when (i) $1\leq k \leq d$ is fixed and (ii) $k=k_d \to \infty$ and $k=\smallO(d)$ as $d \to \infty$.
%, and (iii) $k=k_d \to \infty$ and $k/d \to c \in (0,1)$ as $d \to \infty$.
Relations (\ref{vector_cond:inf}) and (\ref{vector_cond:inf2})
divide the set of all possible values of $(k_d,\mu_d)$ into the region where exact variable selection is possible and
the region where this type of selection is impossible. These relations
describe part of the so-called “phase diagram” which, in the context of variable selection, was
for the first time obtained in \cite{Genovese-2012}, see also \cite{BSt-2017}, \cite{GS-2020}, \cite{Ji-2012}, \cite{MSt-2023}. 
In this work, we establish the analogues of (\ref{vector_cond:inf}) and (\ref{vector_cond:inf2}) for the functional model
at hand, and propose an asymptotically minimax procedure that, when possible, achieves exact recovery of the sparsity pattern of $f$
in model (\ref{model1}) to (\ref{sparsitycond}). The related problem of almost full recovery
of $f$, as studied in \cite{BSt-2017} for $k=1$, will be treated elsewhere.

In the Gaussian white noise model, the problem of selecting active (nonzero) components
of a~sparse additive $d$-variate regression function $f$ from a Sobolev ball was studied
for the univariate components case only, see \cite{BSt-2017} and \cite{ISt-2014}. In this article, assuming that $d$ tends to infinity,
we examine the problem of identifying  active $k$-variate components of a sparse smooth function $f$
when $k$ is an arbitrary integer between 1 and $d$. We first consider the case of fixed $k$ and
then extend the obtained results for the case of growing $k$.

\subsection{Regularity conditions}\label{M}
Consider the sparse Gaussian white noise model given by (\ref{model1}) to (\ref{sparsitycond}).
To ensure that the model at hand is mathematically valid and the problem of interest has a good solution,
some regularity constraints on $f(\tb)=\sum_{u\in {\cal U}_{k,d}} \eta_u f_u(\tb_u)$
need to be imposed. Regression functions in nonparametric regression models are commonly assumed to belong to Sobolev classes of functions.

For $u=\{j_1,\ldots,j_k\}\in{\cal U}_{k,d}$, where $1\leq j_1<\ldots<j_k\leq d,$  define the following sets of vectors:
\begin{gather*}
	\mathbb{Z}_u=\{\lb=(l_1,\ldots,l_d)\in \mathbb{Z}^d: l_j=0\;\mbox{for}\;j\notin u \}, \;
	\mathring{\mathbb{Z}}_u=\{\lb\in \mathbb{Z}_u: l_j\neq 0\;\mbox{for}\;j\in u\}.
\end{gather*}
We also set $\mathring{\mathbb{Z}}_\emptyset=\underbrace{(0,\ldots,0)}_d$, $\mathring{\mathbb{Z}}=\mathbb{Z}\setminus\{0\},$ and
${\mathring{\mathbb{Z}}}^k=\underbrace{\mathring{\mathbb{Z}}\times\ldots\times \mathring{\mathbb{Z}}}_k.$
Then
$\mathbb{Z}^d=\left( \mathring{\mathbb{Z}}\cup\{0\} \right)^d=\bigcup\limits_{u\subseteq\{1,\ldots,d\}}\mathring{\mathbb{Z}}_u.$
Consider the Fourier basis $\{\phi_{\lb}(\tb)\}_{\lb\in{\mathbb{Z}}^d}$ of $L_2^d$  defined as follows:
\begin{equation}
	\begin{split}
		\phi_{\lb}(\tb)&=\prod_{j=1}^d \phi_{l_j}(t_{j}),\quad \lb\in {\mathbb{Z}}^d,
		\\
		\phi_0(t)=1,\quad \phi_l(t)&=\sqrt{2}\cos(2\pi l t),\quad \phi_{-l}(t)=\sqrt{2}\sin(2\pi l t),\quad l>0.
		\label{def:ON_system}
	\end{split}
\end{equation}
Observe that
$\phi_{\lb}(\tb)=\phi_{\lb}(\tb_u)$ for $\lb\in\mathring{\mathbb{Z}}_u$. 
Let function $f$ be determined by (\ref{f}) and satisfy (\ref{orthcon}), and
let $\theta_{\lb}(u)=( f_u,\phi_{\lb})_{L_2^d}$ for $\lb\in\mathring{\mathbb{Z}}_u$. Then
the term $f_u$ can be expressed as
\begin{gather*}
	f_u(\tb_u)=\sum_{\lb\in \mathring{\mathbb{Z}}_u}\theta_{\lb}(u) \phi_{\lb}(\tb_{u})
\end{gather*}
and the entire regression function $f(\tb)=\sum_{u\in {\cal U}_{k,d}} \eta_u f_u(\tb_u)$ takes the form
\begin{gather}\label{f-decomp}
	f(\tb)=\sum_{u\in {\cal U}_{k,d}}\eta_u \sum_{\lb\in \mathring{\mathbb{Z}}_u}\theta_{\lb}(u) \phi_{\lb}(\tb_{u}). 
\end{gather}
For use later on, observe that
only those Fourier coefficients of $f$ that correspond to the orthogonal components $f_u$ are nonzero and
that $\|f_u\|^2_{2}=\sum_{\lb\in \mathring{\mathbb{Z}}_u} \theta_{\lb}^2 (u)$, where $\| \cdot \|_2$ stands for either the norm in $L_2^d$ or $L_2^k$ (the particular choice is clear from the context).

For $u\in{\cal U}_{k,d}$ as above, assume that $f_u(\tb_u)$ belongs to the Sobolev class of $k$-variate functions with integer smoothness parameter $\sigma>0$ for which
the semi-norm $\|\cdot\|_{\sigma,2}$ is defined by
$$\|f_u\|_{\sigma,2}^2=\sum_{i_1=1}^k\ldots\sum_{i_{\sigma}=1}^k \left\|\frac{\partial^{\sigma} f_u  }{\partial t_{j_{i_1}}\ldots\partial t_{j_{i_{\sigma}}  }}   \right\|_{2}^2.$$
Assume, in addition, that $f_u(\tb_u)$ admits 1-periodic $[\sigma]$-smooth extension in each argument to $\mathbb{R}^k$,
i.e., for all derivatives $f_u^{(n)}$ of  order $0\leq n\leq \sigma$, where $f_u^{(0)}=f_u$,
\begin{gather*}
	f_u^{(n)}(t_{j_1},\ldots,t_{j_{i-1}},0,t_{j_{i+1}},\ldots,t_{j_k})=f_u^{(n)}(t_{j_1},\ldots,t_{j_{i-1}},1,t_{j_{i+1}},\ldots,t_{j_k}),
\end{gather*}
for $2\leq i\leq k-1$ with obvious extension for $i=1,k$. Then, the semi-norm $\|f_u\|_{\sigma,2}$ can be defined
for any $\sigma>0$  in terms of the Fourier coefficients $\theta_{\lb}(u)$ as follows:
\begin{gather}
	\|f_u\|_{\sigma,2}^2=\sum_{\lb\in \mathring{\mathbb{Z}}_u } \theta_{\lb}^2(u) c_{\lb}^2,\quad c_{\lb}^2=\left(\sum_{j=1}^d (2\pi l_{j})^2  \right)^{\sigma}
	=\left(\sum_{i=1}^k (2\pi l_{j_i})^2  \right)^{\sigma}.
	\label{def:c_lb}
\end{gather}
Denote by ${\cal F}_{\cb_u}$ the Sobolev ball  of radius 1 with coefficients ${\cb}_u=({c}_{\lb})_{\lb\in\mathring{\mathbb{Z}}_u }$, that is,
\begin{gather*}
	{\cal F}_{\cb_u}=\left \{  f_u(\tb_u)= \sum_{\lb\in\mathring{\mathbb{Z}}_u } \theta_{\lb}(u)\phi_{\lb}(\tb_u), \; \tb_u \in [0,1]^k:
	\sum_{\lb\in \mathring{\mathbb{Z}}_u } \theta_{\lb}^2(u) c_{\lb}^2 \leq 1 \right \},
\end{gather*}
and impose the following component-wise regularity constraints on $f(\tb)=\sum_{u\in {\cal U}_{k,d}} \eta_u f_u(\tb_u)$:
$$f_u\in {\cal F}_{\cb_u}\;\mbox{for all}\; u\in {\cal U}_{k,d}.$$

As seen from {\rm (\ref{def:c_lb})}, the collections $\cb_u=(c_{\lb})_{\lb\in\mathring{\mathbb{Z}}_u}$, $u\in {\cal U}_{k,d}$,
are invariant with respect to a particular choice of the elements of $u\in{\cal U}_{k,d}$ for all $1\leq k\leq d$. Hence, if $u=\{j_1,j_2,\ldots,j_k\}$, $1\leq j_1<\ldots<j_k\leq d$, and $\lb\in \mathring{\mathbb{Z}}_u$, we can write $c_{\lb}=c_{l_{j_1},\ldots,l_{j_k},0,\ldots,0}$.
Thinking of a transformation from $\mathring{\mathbb{Z}}^k$ to $\mathring{\mathbb{Z}}_u$ that maps
$\lb^{\prime}=(l^{\prime}_{1},\ldots,l^{\prime}_{k})\in \mathring{\mathbb{Z}}^k$ to
$\lb=(l_{j_1},\ldots,l_{j_k},0,\ldots,0)\in \mathring{\mathbb{Z}}_u$ according to the rule
$l_j = 0$ for $ j \notin u$ and $l_{j_s}=l^{\prime}_s$ for $s=1,\ldots,k,$
we may represent each collection  $\cb_u=(c_{\lb})_{\lb\in\mathring{\mathbb{Z}}_u}$
as $\cb_k=(c_{\lb})_{\lb\in\mathring{\mathbb{Z}}^k}$.
Clearly, the collections of functions ${\cal F}_{\cb_u}$ are isomorphic for all $u$ of cardinality $k$.
This observation will be used to state and prove the main results.

\subsection{Problem statement}\label{PS}
To ensure that the problem of identifying the
active components of $f$ in model (\ref{model1}) to (\ref{sparsitycond}) is meaningful, we would naturally require that
the components $f_u$ of $f$ are not too ``small''.
For this purpose, following \cite{BSt-2017} and \cite{ISt-2014},
for $u\in {\cal U}_{k,d}$, $1\leq k\leq d$, and $r>0$,
define
\begin{gather*}
	{\mathring{\cal F}}_{\cb_u}(r)=\left\{f_u\in {\cal F}_{\cb_u} : {\|f_u\|}_2\geq r\right\}
\end{gather*}
and consider testing a simple null hypothesis versus a family of alternatives: 
\begin{gather}\label{test1}
	H_{0,u}:f_u=0\quad\mbox{vs.}\quad H^\varepsilon_{1,u}: f_u\in \mathring{{\cal F}}_{\cb_u}(r_{\varepsilon,k})
\end{gather}
for some positive family $r_{\varepsilon,k}\to 0$ as $\varepsilon\to 0$. 
The hypothesis testing problem in~(\ref{test1}), known in the literature as the signal detection problem,
has been studied in \cite{IS-2015}. In the present context of sparse signal recovery,
this is an auxiliary problem that makes it possible to obtain the conditions when exact selection of active components of $f$ is possible and
when this type of selection is impossible. Additionally, we shall use the
asymptotically minimax test statistics from the above signal detection problem (see Theorem 3.1 of \cite{IS-2015}) to design an exact selector.

Recall that two hypotheses $H_{0}$ and $H_{1}$ separate asymptotically if there exists a consistent test procedure for testing
$H_{0}$ against $H_{1}.$ When $H_{0,u}$ and $H^\varepsilon_{1,u}$ as in (\ref{test1}) separate asymptotically, we say that $f_u$ is \textit{detectable}.
It will be seen that only detectable components $f_u$ of $f$ can be correctly identified. 
Moreover, the quantity that is crucial for establishing the \textit{sharp selection boundary} (see the quantity $a_{\varepsilon,u} (r_{\varepsilon,k})$ in inequalities (\ref{cond:inf}) and (\ref{cond:inf2}) in Section \ref{FixedK}),
turns out to be exactly
the quantity that defines the sharp detection boundary in the signal detection problem (\ref{test1})
(see inequalities (\ref{testing}) in Section \ref{Selector}).

For a positive family $r_{\varepsilon,k}\to 0$ as $\varepsilon\to 0$, define the class of sparse multivariate functions of our interest by
\begin{multline*}
	{\cal F}^d_{k,\sigma}(r_{\varepsilon,k})=\left\{f: f(\tb)=\sum_{u\in {\cal U}_{k,d}} \eta_u f_u(\tb_u), f_u\;\mbox{satisfies}\;(\ref{orthcon}),\right.\\
	\left.f_u\in \mathring{\cal F}_{\cb_u}(r_{\varepsilon,k}), u\in {\cal U}_{k,d}, \etab=(\eta_u)_{u\in {\cal U}_{k,d}}\in \mathcal{H}^d_{k,\beta} \right\},
\end{multline*}
where the dependence of ${\cal F}^d_{k,\sigma}(r_{\varepsilon,k})$ on the smoothness parameter $\sigma$ is hidden in the coefficients $\cb_u=(c_{\lb})_{\lb\in\mathring{\mathbb{Z}}_u}$.
In this work, we obtain the sharp selection boundary that allows us to separate identifiable components $f_u$ of
a signal $f\in {\cal F}^d_{k,\sigma}(r_{\varepsilon,k})$ in model (\ref{model1}) from non-identifiable ones, that is, we establish the conditions on
$r_{\varepsilon,k}$ for the possibility of exact identification of the nonzero components $f_u$ of~$f$.
Next, we construct
an adaptive (free of $\beta$) estimator $\hat{\etab}=\hat{\etab}(X_\varepsilon)\in\{0,1\}^{d\choose k}$ of $\etab$ attaining this boundary with the property
\begin{gather}\label{r1}
	\limsup_{\varepsilon\to 0}\sup_{\etab\in \mathcal{H}^d_{k,\beta}}\sup_{f\in  {\cal F}^d_{k,\sigma}(r_{\varepsilon,k})} {\rm E}_{f,\etab}|\hat{\etab}-\etab| =0,
\end{gather}
which holds true for all $\beta\in(0,1)$ and all $\sigma>0$, where $${\rm E}_{f,\etab}|\hat{\etab}-\etab| ={\rm E}_{f,\etab}\left(\sum_{u\in{\cal U}_{k,d}}|\hat{\eta}_u-\eta_u|\right)$$ is the \textit{Hamming risk} of $\hat{\etab}$ for given $f$ and $\etab$. Finally, we show that for
those values of $r_{\varepsilon,k}$ that fall below the selection boundary,
one has
\begin{gather}\label{r2}
	\liminf_{\varepsilon\to 0}\inf_{\tilde{\etab}}\sup_{\etab\in\mathcal{H}^d_{k,\beta}}\sup_{f\in  {\cal F}^d_{k,\sigma}(r_{\varepsilon,k})} {\rm E}_{f,\etab}|\tilde{\etab}-\etab|>0,
\end{gather}
that is, exact recovery of the signal $f\in {\cal F}^d_{k,\sigma}(r_{\varepsilon,k})$ in model (\ref{model1}) to (\ref{sparsitycond}) is impossible.

The initial model given by (\ref{model1}) to (\ref{sparsitycond}) can be restated in terms of the (nonzero) Fourier coefficients
of $f(\tb)=\sum_{u\in {\cal U}_{k,d}}\eta_u f_u(\tb_u)$ as follows. Note that the orthonormal system (\ref{def:ON_system}) satisfies $\{\phi_{\lb}(\tb)\}_{\lb \in  \mathbb{Z}^d} = \bigcup\limits_{u\subseteq\{1,\ldots,d\}} \{\phi_{\lb}(\tb_{u})\}_{\lb \in  \mathring{\mathbb{Z}}_u}$. Next, recall the orthogonality conditions (\ref{orthcon}) and let $X_{\lb}=X_\varepsilon(\phi_{\lb})$ be the $\lb$th empirical Fourier coefficient, that is,
\begin{align*}
	X_{\lb}&=\int_{[0,1]^d}\phi_{\lb}(\tb)\, dX_{\varepsilon}(\tb)=
	\sum_{v\in {\cal U}_{k,d}} \eta_v \int_{[0,1]^d} f_v(\tb_v)\phi_{\lb}(\tb)\, d\tb+\varepsilon\int_{[0,1]^d}\phi_{\lb}(\tb)\, dW(\tb)\\
	&=\eta_u\theta_{\lb}(u)+\varepsilon\xi_{\lb},\quad \lb \in  \mathring{\mathbb{Z}}_u,
\end{align*}
where $\theta_{\lb}(u)= (f_u,\phi_{\lb})_{L_2^d}$ is
the $\lb$th Fourier coefficient of $f_u$, and random variables $\xi_{\lb}=$ $\int_{[0,1]^d}\phi_{\lb}(\tb)\, dW(\tb)$
are iid standard normal for ${\lb \in  \mathring{\mathbb{Z}}_u}$ and $u\in {\cal U}_{k,d}$. Then,
the initial model described by (\ref{model1}) to (\ref{sparsitycond}) turns to the sequence space model
\begin{gather}\label{model2}
	X_{\lb}=\eta_u\theta_\lb(u)+\varepsilon \xi_{\lb},\quad \lb\in \mathring{\mathbb{Z}}_u, \quad u\in {\cal U}_{k,d},
\end{gather}
where $\theta_{\lb}(u)$ and $\xi_{\lb}$ are as above,
and $\etab=(\eta_u)_{u\in {\cal U}_{k,d}}\in {\cal H}^d_{k,\beta}$.
The problem of testing $H_{0,u}$ against $H^{\varepsilon}_{1,u}$ in (\ref{test1}) is therefore equivalent to that of testing
\begin{gather*}
	\mathbf{H}_{0,u}: \thetab_u=\boldsymbol{0} \quad\mbox{vs.}\quad \mathbf{H}^{\varepsilon}_{1,u}: \thetab_u \in \mathring{\Theta}_{\cb_u}(r_{\varepsilon,k}),
\end{gather*}
where $\thetab_u=(\theta_{\lb}(u))_{\lb\in \mathring{\mathbb{Z}}_u}$ and
\begin{equation}
	\mathring{\Theta}_{\cb_u}(r_{\varepsilon,k})=\left\{\thetab_u=(\theta_{\lb}(u))_{\lb\in \mathring{\mathbb{Z}}_u}\in l_2(\mathbb{Z}^d):
	\sum_{\lb \in \mathring{\mathbb{Z}}_u}\theta_{\lb}^2(u)c_{\lb}^2\leq 1, \sum_{\lb \in \mathring{\mathbb{Z}}_u}\theta_{\lb}^2(u)\geq r_{\varepsilon,k}^2\right\}.
	\label{def:Theta_ring}
\end{equation}
Observe that $\mathring{\Theta}_{\cb_u}(r_{\varepsilon,k})=\emptyset$ when $r_{\varepsilon,k}>1/c_{\varepsilon,0},$ where,
recalling (\ref{def:c_lb}), $c_{\varepsilon,0}:=\inf_{\lb\in \mathring{\mathbb{Z}}_u}c_{\lb}=(2\pi)^{\sigma}k^{\sigma/2}$.
Therefore, in what follows, we are interested in the case when $r_{\varepsilon,k}\in (0,(2\pi)^{-\sigma}k^{-\sigma/2}).$

The recovery of $\etab=(\eta_u)_{u\in {\cal U}_{k,d}}$ in the sequence space model (\ref{model2}) will be named the \textit{variable selection problem}.
In terms of model (\ref{model2}),
relations (\ref{r1}) and (\ref{r2}) take the form
\begin{gather}\label{r11}
	\limsup_{\varepsilon\to 0}\sup_{\etab\in\mathcal{H}^d_{k,\beta}}\sup_{\thetab\in  \Theta^d_{k,\sigma}(r_{\varepsilon,k})} {\rm E}_{\thetab,\etab}|\hat{\etab}-\etab| =0,
\end{gather}
and
\begin{gather}\label{r22}
	\liminf_{\varepsilon\to 0}\inf_{\tilde{\etab}}\sup_{\etab\in\mathcal{H}^d_{k,\beta}}\sup_{\thetab\in  {\Theta}^d_{k,\sigma}(r_{\varepsilon,k})} {\rm E}_{\thetab,\etab}|\tilde{\etab}-\etab|>0.
\end{gather}
Here,  $\hat{\etab}$ and $\tilde{\etab}$ are estimators of $\etab$ based on  $\{\Xb_{\!u}\}_{u\in {\cal U}_{k,d}}$, where
$\Xb_{\!u}=(X_{\lb})_{\lb\in \mathring{\mathbb{Z}}_u}$, and
\begin{gather*}
	{\Theta}^d_{k,\sigma}(r_{\varepsilon,k})=\left\{\thetab: \thetab=(\thetab_{u})_{u\in {\cal U}_{k,d}},\;\mbox{where}\;
	\thetab_u=(\theta_{\lb}(u))_{\lb\in \mathring{\mathbb{Z}}_u} \in \mathring{\Theta}_{\cb_u}(r_{\varepsilon,k})  \right\}.  %\label{Theta33}
\end{gather*}
Below, an estimator $\hat{\etab}\in\{0,1\}^{d\choose k}$ of $\etab$ based on  $\{\Xb_{\!u}\}_{u\in {\cal U}_{k,d}}$ will be called a~\textit{selector}.
An \textit{exact selector}~$\hat{\etab}$ will be defined as a selector satisfying (\ref{r11}).
The limiting relations (\ref{r11}) and (\ref{r22}) will be referred to as the \textit{upper bound} on the maximum Hamming risk of $\hat{\etab}$
and the \textit{lower bound} on the minimax Hamming risk, respectively.

\subsection{Notation}\label{N}
The majority of limits in the article are taken as $\varepsilon\to 0$.
We use the notation $a_\varepsilon \asymp b_\varepsilon$ when
$0 < \liminf_{\varepsilon \to 0} {a_\varepsilon}/{b_\varepsilon} \leq \limsup_{\varepsilon\to 0} {a_\varepsilon}/{b_\varepsilon} < \infty$
and the notation $a_\varepsilon \sim b_\varepsilon$ when
$\lim_{\varepsilon \to 0} {a_\varepsilon}/{b_\varepsilon} =1.$
For positive $a_\varepsilon$ and $b_\varepsilon$, we shall occasionally write $a_\varepsilon\ll b_\varepsilon$
when $a_\varepsilon=\smallO(b_\varepsilon)$ as $\varepsilon\to 0$.
For a real number $x$,  $x_+ = \max (0,x)$.
We use the symbol $\log x$ for the natural (base $e$) logarithm of the number $x$.
A~noncentral chi-square random variable with $\nu$ degrees of freedom and
noncentrality parameter $\lambda$ will be denoted by $\chi^2_{\nu}(\lambda)$.

\section{Construction of an exact selector}\label{Selector}
For $u\in {\cal U}_{k,d}$ and  $r_{\varepsilon,k}>0$, consider the quantity
\begin{equation}
	a^2_{\varepsilon,u} (r_{\varepsilon,k})= \frac{1}{2 \varepsilon^4} \inf_{\thetab_u \in \mathring{\Theta}_{\cb_u}(r_{\varepsilon,k})} \sum_{\lb\in \mathring{\mathbb{Z}}_u} \theta^4_\lb(u),
	\label{def:a2}
\end{equation}
which plays a key role in the minimax theory of hypothesis testing.
In particular, the problem of testing $\mathbf{H}_{0,u}:\thetab_u=\boldsymbol{0}$ versus $\mathbf{H}^{\varepsilon}_{1,u}:\thetab_u \in \mathring{\Theta}_{\cb_u}(r_{\varepsilon,k})$ is meaningful when
$r_{\varepsilon,k}/r^{\prime}_{\varepsilon,k}\to \infty$, where the \textit{detection boundary} $r^{\prime}_{\varepsilon,k}$ is determined by $a_{\varepsilon,u}(r^{\prime}_{\varepsilon,k})\asymp \sqrt{\log{d\choose k}  }$.
If the hypotheses $\mathbf{H}_{0,u}$ and $\mathbf{H}^{\varepsilon}_{1,u}$ separate asymptotically when
$\liminf_{\varepsilon\to 0}r_{\varepsilon,k}/r^{\prime}_{\varepsilon,k}>1,$
and merge asymptotically when
$\limsup_{\varepsilon\to 0}r_{\varepsilon,k}/r^{\prime}_{\varepsilon,k}<1,$
then the detection boundary $r^{\prime}_{\varepsilon,k}$ is said to be \textit{sharp}.
The \textit{sharp detection boundary} in the above hypothesis testing problem is determined in terms of $a^2_{\varepsilon,u} (r_{\varepsilon,k})$ by
the  inequalities (see Theorem 3.1 of~\cite{IS-2015})
\begin{gather}\label{testing}
	{\rm (a)} \;\liminf_{\varepsilon \to 0}  \frac{a_{\varepsilon,u}(r_{\varepsilon,k})}{\sqrt{\log {d\choose k}}} > \sqrt{2}
	\quad \mbox{and}\quad {\rm (b)}\; \limsup_{\varepsilon \to 0}  \frac{a_{\varepsilon,u}(r_{\varepsilon,k})}{\sqrt{\log {d\choose k}}} < \sqrt{2}.
\end{gather}
Specifically, when inequality (a) holds,
the hypotheses $\mathbf{H}_{0,u}$ and $\mathbf{H}^{\varepsilon}_{1,u}$ separate asymptotically,
whereas when inequality (b)  holds, $\mathbf{H}_{0,u}$ and $\mathbf{H}^{\varepsilon}_{1,u}$ merge asymptotically (i.e., there is no consistent test procedure
for testing $\mathbf{H}_{0,u}$ versus $\mathbf{H}^{\varepsilon}_{1,u}$).
In the present context of variable selection, the quantity $a_{\varepsilon,u}(r_{\varepsilon,k})$ is also of great importance.

Observe that $a_{\varepsilon,u} (r_{\varepsilon,k})$ is a nondecreasing function of its argument that possesses a kind of ``continuity'' property. Namely, for all $\gamma >0$, there exist $\varepsilon^*>0$ and $\delta^* >0$ such that (see Section 5.2.3 of~\cite{NGF-2003})
\begin{equation}
	a_{\varepsilon,u} (r_{\varepsilon,k}) \leq a_{\varepsilon,u} ((1+\delta)r_{\varepsilon,k}) \leq (1+\gamma) a_{\varepsilon,u} (r_{\varepsilon,k}), \quad \forall\, \varepsilon \in (0,\varepsilon^*), \forall\, \delta \in (0,\delta^*).
	\label{def:continuity_of_a}
\end{equation}
These and other general facts of the minimax hypothesis testing theory can be found in
a series of review articles \cite{I1993a}--\cite{I1993c} and monograph \cite{NGF-2003}.

Suppressing for brevity dependence on $u$, denote the minimizing sequence in (\ref{def:a2}) by $(\theta^*_\lb (r_{\varepsilon,k}))_{\lb \in \mathring{\mathbb{Z}}_u}$, that is,
\begin{equation}
	a^2_{\varepsilon,u} (r_{\varepsilon,k})= \frac{1}{2 \varepsilon^4} \sum_{\lb\in \mathring{\mathbb{Z}}_u} \left[\theta^{*}_\lb (r_{\varepsilon,k})\right]^4
	\label{def:a2_with_theta_star}
\end{equation}
and let $r^*_{\varepsilon,k}>0$ be determined by, cf. formula (14) in \cite{ISt-2014},
\begin{equation}
	a_{\varepsilon,u} (r^*_{\varepsilon,k}) = (1+\sqrt{1-\beta}) \sqrt{2 \log {d\choose k}}.
	\label{def:r_star}
\end{equation}
For the purpose of constructing an exact selector $\hat{\etab}$ satisfying (\ref{r11}), consider the weighted $\chi^2$-type statistics 
\begin{equation}
	S_{u} = \sum_{\lb \in \mathring{\mathbb{Z}}_u} \omega_\lb (r^*_{\varepsilon,k}) \left[ \left( {X_\lb}/{\varepsilon}  \right)^2 -1 \right], \quad u \in {\cal U}_{k,d},
	\label{def:S_u}
\end{equation}
where
\begin{equation}
	\omega_\lb (r_{\varepsilon,k}) = \frac{1}{2 \varepsilon^2} \frac{[\theta^*_\lb (r_{\varepsilon,k})]^2}{a_{\varepsilon,u}(r_{\varepsilon,k})}, \quad \lb \in \mathring{\mathbb{Z}}_u.
	\label{def:omega}
\end{equation}
Due to (\ref{def:a2_with_theta_star}), one has
\begin{equation}
	\sum_{\lb \in \mathring{\mathbb{Z}}_u} \omega_\lb^2 (r_{\varepsilon,k}) ={1}/{2} \quad \mbox{for all}\;\; r_{\varepsilon,k} > 0.
	\label{eq:sum_omega_is_half}
\end{equation}
The statistics $S_{u}$ in (\ref{def:S_u}), which appeared earlier in \cite{IS-2015} in connection with the signal detection problem~(\ref{test1}),
depend on the sparsity index $\beta$ through the weights $\omega_\lb (r^*_{\varepsilon,k}),$ $\lb\in \mathring{\mathbb{Z}}_u$.

It is known (see, for example, Remark 2.1 of \cite{IS-2015} or Theorem 4 of \cite{IS-2005}) that
\begin{equation}
	[\theta^*_\lb (r_{\varepsilon,k})]^2
	=a_0^2 \left( 1-(c_{\lb}/T)^2 \right)_+,
	\label{def:theta_star}	
\end{equation}
where the collection $(c_\lb)_{{\lb \in \mathring{\mathbb{Z}}_u}}$ is defined in (\ref{def:c_lb}) and the quantities $a_0^2=a_{0,\varepsilon,k}^2$ and $T=T_{\varepsilon,k}$ are determined by
\begin{equation*}
	a_0^2\sum_{\lb \in \mathring{\mathbb{Z}}_u}\left(1-(c_{\lb}/T)^2  \right)=r_{\varepsilon,k}^2,\quad
	a_0^2\sum_{\lb \in \mathring{\mathbb{Z}}_u}c_{\lb}^2\left(1-(c_{\lb}/T)^2  \right)=1,\quad T\geq r^{-1}_{\varepsilon,k}\to \infty.
\end{equation*}
For fixed $k$, by using Theorem 4 of \cite{IS-2005} and its proof,
we can easily get the asymptotic expressions for $a_0^2$ and $T$ on the right-hand side of (\ref{def:theta_star}), and then, for $\lb\in \mathring{\mathbb{Z}}_u$, by recalling (\ref{def:c_lb}), we obtain
as $\varepsilon\to 0$
\begin{equation}\label{thetal}
	[\theta^*_\lb (r_{\varepsilon,k})]^2 \sim
	\frac{r_{\varepsilon,k}^{2+{k}/{\sigma}} 2^k \pi^{{k}/{2}}(k+2\sigma)\Gamma\left(1+{k}/{2}\right) }{2\sigma\left(1+4\sigma/k  \right)^{{k}/{(2\sigma)}}}  \left(1-\left(\sum\nolimits_{j=1}^d (2\pi l_j)^2  \right)^{\sigma}\frac{r^2_{\varepsilon,k} }{(1+4\sigma/k)}  \right)_+.
\end{equation}	
Since
\begin{equation*}
\#\{\lb \in \mathring{\mathbb{Z}}_u: \theta^*_\lb (r_{\varepsilon,k})\neq 0\}
=\#\left\{\lb \in \mathring{\mathbb{Z}}_u:\left(\sum_{j=1}^d l_j^2  \right)^{1/2}< (1+4\sigma/k)^{1/(2\sigma)}/(2\pi r_{\varepsilon,k}^{1/\sigma})\right\}
\end{equation*}
is of the same order of magnitude as the integer-volume of the $l_2$-ball of radius $(1+4\sigma/k)^{1/(2\sigma)}/(2\pi r_{\varepsilon,k}^{1/\sigma})$ in $\mathbb{R}^k$,
which is  $\mathcal{O}(r_{\varepsilon,k}^{-{k}/{\sigma}})$, it follows that
\begin{gather*}
	\#\{\lb \in \mathring{\mathbb{Z}}_u: \theta^*_\lb (r_{\varepsilon,k})\neq 0\}\asymp r_{\varepsilon,k}^{-{k}/{\sigma}},\quad \varepsilon\to 0,
\end{gather*}
when both $k$ is fixed and $k=k_\varepsilon\to \infty$.

We can modify the asymptotic expression in (\ref{thetal}) for $k=k_\varepsilon\to \infty$.
Noting that $\Gamma(x+1)\sim \sqrt{2\pi}x^{x+1/2}e^{-x}$ and $(1+1/x)^x\sim e$ as $x\to \infty$ gives
\begin{gather}\label{theta2}
	[\theta^*_\lb (r_{\varepsilon,k})]^2\sim
	\frac{r_{\varepsilon,k}^{2+{k}/{\sigma}  }\pi^{1/2}(2\pi k/e)^{k/2} k^{3/2} }{2  \sigma e^2}  \left(1-\left(\sum\nolimits_{j=1}^d (2\pi l_j)^2  \right)^{\sigma}\frac{r^2_{\varepsilon,k} }{(1+4\sigma/k)}  \right)_+
\end{gather}
as $\varepsilon\to 0$. In both cases, when $k$ is fixed and $k=k_\varepsilon\to \infty$, in view of (\ref{def:S_u}) and~(\ref{def:omega}),  every statistic $S_u$, $u\in {\cal U}_{k,d}$,
consists of $\mathcal{O}(r_{\varepsilon,k}^{-{k}/{\sigma}})$ nonzero terms.

It is also known that for fixed $k$, the sharp asymptotics of $a_{\varepsilon,u}(r_{\varepsilon,k})$, $u\in {\cal U}_{k,d}$, as $\varepsilon\to 0$  are given by (see Theorem 4 of \cite{IS-2005})
\begin{gather}\label{aek}
	a_{\varepsilon,u}(r_{\varepsilon,k})\sim C(\sigma,k)r_{\varepsilon,k}^{2+{k}/{(2\sigma)}}\varepsilon^{-2},\quad
	C^2(\sigma,k)=\frac{\pi^k (1+{2\sigma}/{k})\Gamma(1+k/2)}{(1+4\sigma/k)^{1+{k}/{(2\sigma)}}\Gamma^k(3/2)}
\end{gather}
and, when $k=k_\varepsilon\to \infty$, cf. formula (39) in \cite{IS-2005},
\begin{gather}\label{aek1}
	a_{\varepsilon,u}(r_{\varepsilon,k})\sim \left({2\pi k}/{e}\right)^{{k}/{4}}e^{-1}(\pi k)^{1/4}r_{\varepsilon,k}^{2+{k}/{(2\sigma)}}\varepsilon^{-2}.
\end{gather}

In case of known $\beta$,
in order to estimate a vector $\etab=(\eta_u)_{u\in {\cal U}_{k,d}}\in{\cal H}^d_{k,\beta}$,
we could use a selector $\boldsymbol{\check{\eta}}=(\check{\eta}_u)_{u\in {\cal U}_{k,d}}$, depending on the data $\{\Xb_{\!u}\}_{u\in {\cal U}_{k,d}}$
through the statistics $S_{u}$, ${u\in {\cal U}_{k,d}}$,
of the form, cf.~(\ref{vector_selector}),
\begin{equation}
	\check{\eta}_u = \ind{ S_{u} > \sqrt{(2+\epsilon) \log {d\choose k} }},\quad u\in {\cal U}_{k,d},
	\label{def:selector_for_beta_known}
\end{equation}
with some $\epsilon=\epsilon_\varepsilon \to 0$ as $\varepsilon\to 0$.
However, the statistics $S_u$ depend on $\beta$, and the sparsity parameter $\beta$ is usually unavailable.
In order to adapt $S_u$ to an unknown value of $\beta$, as in Section 3.1 of \cite{ISt-2014},  consider a grid of equidistant points
\begin{gather*}
	\beta_{k,m} = m \rho_k,\quad m=1,\ldots,M_k,
\end{gather*}
where $\rho_k>0$ is a small number that depends on $\varepsilon$, $M_k = \lfloor 1/\rho_k \rfloor$, and assume that $\log M_k = \smallO \left( \log {d\choose k} \right)$ as
$\varepsilon \to 0$.
Next, for $m=1,\ldots,M_k$, define $r^*_{\varepsilon,k,m} >0$ by
\begin{equation}
	a_{\varepsilon,u} (r^*_{\varepsilon,k,m}) = (1+\sqrt{1-\beta_{k,m}}) \sqrt{2 \log {d\choose k} }
	\label{def:r_star_ekm}
\end{equation}
and for every $u \in {\cal U}_{k,d}$ consider the statistics
\begin{equation}
	S_{u,m} = \sum_{\lb \in \mathring{\mathbb{Z}}_u} \omega_\lb (r^*_{\varepsilon,k,m}) \left[ \left( {X_\lb}/{\varepsilon}  \right)^2 -1 \right], \quad m=1,\ldots,M_k.
	\label{def:S_um}
\end{equation}
Now, we propose a selector $\boldsymbol{\hat{\eta}}=(\hat{\eta}_u)_{u \in {\cal U}_{k,d} }$ with components, cf. (\ref{def:selector_for_beta_known}),
\begin{equation}
	\hat{\eta}_u = \max_{1 \leq m \leq M_k} \hat{\eta}_{u,m}, \quad  \hat{\eta}_{u,m} = \ind{ S_{u,m} > \sqrt{(2+\epsilon) \left( \log {d\choose k} + \log M_k \right)} },
	\label{def:selector_for_beta_unknown}
\end{equation}
where $\epsilon>0$ depends on $\varepsilon$ and satisfies
\begin{gather}\label{delta}
	\epsilon\to 0\quad\mbox{and}\quad \epsilon\log {d\choose k}\to \infty,\quad \mbox{as}\;\; \varepsilon\to 0.
\end{gather}
In this work, we shall deal with the cases when $d=d_\varepsilon\to \infty$ and $k$ is either fixed or $k=k_\varepsilon\to \infty$ and $k=\smallO(d)$ as $\varepsilon\to \infty$.

Let us explain the idea behind the selector $\hat{\etab}=(\hat{\eta}_u)_{u \in {\cal U}_{k,d} }$
defined by (\ref{def:selector_for_beta_unknown}).
For $u\in {\cal U}_{k,d}$,  $\hat{\eta}_u$~identifies the component $f_u$  as active, when at least one
of the statistics $ S_{u,m}$, $m=1,\ldots,M_k$, detects it. Therefore,  the probability of falsely non-identifying
$f_u$ as active by means of $\hat{\eta}_u$ does not exceed the probability $P_{\eta_u=1}(\hat{\eta}_{u,m_0} = 0)$
for some $\beta_{k,m_0}$ close to the true (but unknown) value of $\beta$.
Also, the probability of falsely identifying $f_u$ as active is bounded from above by the sum
$\sum_{m=1}^{M_k} P_{\eta_u=0}(\hat{\eta}_{u,m} =1)$, which is small by the choice of the threshold
(for more details, see the proof of  Theorem \ref{th:UB_fixedK}).
It should be understood that the proposed selector $\hat{\etab}$ depends on $\sigma>0$, and therefore, we have a whole class of adaptive selectors indexed by $\sigma$.

\begin{remark}\label{remark2}
	Consider the ellipsoid with a small central ball removed, $\mathring{\Theta}_{\cb_u}(r_{\varepsilon,k})$, defined in {\rm (\ref{def:Theta_ring})}.
	When $\lb$ belongs to $\mathring{\mathbb{Z}}^k$ instead of $\mathring{\mathbb{Z}}_u$,
	we shall use the notation	$\mathring{\Theta}_{\cb_k}(r_{\varepsilon,k})$, that is,
	\begin{equation*}
		\mathring{\Theta}_{\cb_k}(r_{\varepsilon,k})=\left\{\thetab_k=(\theta_{\lb})_{\lb\in \mathring{\mathbb{Z}}^k}\in l_2(\mathbb{Z}^d):
		\sum_{\lb \in \mathring{\mathbb{Z}}^k}\theta_{\lb}^2 c_{\lb}^2\leq 1, \sum_{\lb \in \mathring{\mathbb{Z}}^k}\theta_{\lb}^2\geq r_{\varepsilon,k}^2\right\}.
	\end{equation*}
	Likewise,
	when $\lb$ belongs to $\mathring{\mathbb{Z}}^k$ instead of $\mathring{\mathbb{Z}}_u$ in {\rm (\ref{def:S_u})} and {\rm (\ref{def:S_um})},
	we shall write $S_{k}$ and $ S_{k,m}$ in place of $S_u$ and $S_{u,m}$, respectively.
	Finally,
	when $\lb$ belongs to $\mathring{\mathbb{Z}}^k$ instead of $\mathring{\mathbb{Z}}_u$ in  {\rm (\ref{def:a2})}, the value $a^2_{\varepsilon,u}(r_{\varepsilon,k})$ of the extremal problem  will be denoted by $a^2_{\varepsilon,k}(r_{\varepsilon,k})$.
	Below, rather than testing $\mathbf{H}_{0,u}:\thetab_u=0$ versus $\mathbf{H}^{\varepsilon}_{1,u}:\thetab_u\in \mathring{\Theta}_{\cb_u}(r_{\varepsilon,k})$, where the alternative hypotheses
	$\mathbf{H}^{\varepsilon}_{1,u}$
	are of the same form for all $u\in{\cal U}_{k,d}$,
	we shall consider a unified problem of testing
	$\mathbf{H}_{0,k}:\thetab_k=0$ versus $\mathbf{H}^{\varepsilon}_{1,k}: \thetab_k\in \mathring{\Theta}_{\cb_k}(r_{\varepsilon,k})$. 
\end{remark}

\section{Main results}\label{MRes}

In this section, we show that, under certain model assumptions, the selector $\hat{\etab}=(\hat{\eta}_u)_{u \in {\cal U}_{k,d} }$ defined by
(\ref{def:selector_for_beta_unknown}) 
achieves exact selection in model {\rm (\ref{model2})}, and that it is the best possible selector in the asymptotically minimax sense.
We consider the case when $k$ is fixed and the case when $k=k_\varepsilon\to \infty$ and $k=\smallO(d)$ as $\varepsilon\to 0$.

\subsection{The case of fixed $k$}\label{FixedK}
Recall the function $a_{\varepsilon,k}(r_{\varepsilon,k})$, $r_{\varepsilon,k}>0$, defined in (\ref{def:a2}):
\begin{equation*}
	a^2_{\varepsilon,k} (r_{\varepsilon,k})= \frac{1}{2 \varepsilon^4} \inf_{\thetab_k \in \mathring{\Theta}_{\cb_k}(r_{\varepsilon,k})} \sum_{\lb\in \mathring{\mathbb{Z}}^k} \theta^4_\lb.
\end{equation*}
The following result describes the conditions that are required for the selector $\hat{\etab}$ in  (\ref{def:selector_for_beta_unknown}) to be exact.

\begin{theorem} \label{th:UB_fixedK}
	Let $k\in\{1,\ldots,d\}$, $\beta\in(0,1)$, and $\sigma>0$ be fixed numbers, and let $d=d_\varepsilon\to \infty$
	and  $\log {d\choose k}=\smallO \left( \log \varepsilon^{-1} \right)$ as $\varepsilon\to 0$.
	Let the quantity
	$r_{\varepsilon,k}>0$ be such that 
	\begin{equation}
		\liminf_{\varepsilon \to 0} \frac{a_{\varepsilon,k}(r_{\varepsilon,k})}{\sqrt{\log {d\choose k}}} > \sqrt{2} (1+\sqrt{1-\beta}).
		\label{cond:inf}
	\end{equation}
	Then the selector $\boldsymbol{\hat{\eta}}=(\hat{\eta}_u)_{u \in {\cal U}_{k,d} }$ given by {\rm (\ref{def:selector_for_beta_unknown})} satisfies
	\begin{equation*}
		\limsup_{\varepsilon\to 0}	\sup_{\boldsymbol{\eta} \in \mathcal{H}^d_{k,\beta}} \sup_{\thetab\in \Theta_{k,\sigma}^d(r_{\varepsilon,k})} \E{\thetab,\etab} |\boldsymbol{\eta} - \boldsymbol{\hat{\eta}}| = 0.
	\end{equation*}
\end{theorem}

Condition (\ref{cond:inf}) imposed on $r_{\varepsilon,k}$ is the identifiability condition; it ensures that
for every $u\in {\cal U}_{k,d}$
the ``signal'' described by the alternative $\mathbf{H}^{\varepsilon}_{1,u}:\thetab_u \in \mathring{\Theta}_{\cb_u}(r_{\varepsilon,k})$
is identifiable (and also detectable).
Theorem \ref{th:UB_fixedK} says that, under condition (\ref{cond:inf}),
the selection procedure based on $\boldsymbol{\hat{\eta}}$ reconstructs all nonzero components of a vector $\etab\in \mathcal{H}^d_{k,\beta}$
and thus provides exact recovery of $(\eta_u\thetab_u)_{u\in {\cal U}_{k,d}}$, uniformly in $\beta\in(0,1),$ $
\sigma>0$, and over $ \mathcal{H}^d_{k,\beta}$ and $\Theta_{k,\sigma}^d(r_{\varepsilon,k})$. In other words,
the proposed selector $\boldsymbol{\hat{\eta}}$ achieves exact selection in model {\rm (\ref{model2})}.
Theorem \ref{th:UB_fixedK} extends Theorem 1 of \cite{ISt-2014} from $k=1$ to the case $1\leq k\leq d$.

The next theorem shows that if the identifiability condition (\ref{cond:inf}) is not satisfied,
the minimax Hamming risk is strictly positive (in the limit),  and thus
exact selection in model (\ref{model2}) is impossible.

\begin{theorem} \label{theorem:LB_figed_k}
	Let $k\in\{1,\ldots,d\}$, $\beta\in(0,1)$, and $\sigma>0$ be fixed numbers, and let $d=d_\varepsilon\to \infty$
	and $\log {d\choose k}=\smallO \left(\log \varepsilon^{-1} \right)$ as $\varepsilon\to 0$.
	Let the quantity $r_{\varepsilon,k}>0$ be such that
	\begin{equation}
		\limsup_{\varepsilon \to 0} \frac{a_{\varepsilon,k}(r_{\varepsilon,k})}{\sqrt{\log {d\choose k}}} < \sqrt{2} (1+\sqrt{1-\beta}).
		\label{cond:inf2}
	\end{equation}
	Then
	\begin{equation*}
		\liminf_{\varepsilon \to 0}\inf_{\tilde{\etab}} \sup_{\boldsymbol{\eta} \in \mathcal{H}^d_{k,\beta}} \sup_{\thetab\in \Theta_{k,\sigma}^d(r_{\varepsilon,k})} \E{\thetab,\etab} |\boldsymbol{\eta} - \boldsymbol{\tilde{\eta}}|>0,
	\end{equation*}
	where the infimum is taken over all selectors $\tilde{\etab} = (\tilde{\eta}_u)_{u \in \mathcal{U}_{k,d}}$ of a vector $\etab$ in model {\rm (\ref{model2})}. 
\end{theorem}

Theorems \ref{th:UB_fixedK} and  \ref{theorem:LB_figed_k}  imply that the selector $\boldsymbol{\hat{\eta}}=(\hat{\eta}_u)_{u \in \mathcal{U}_{k,d}}$
in {\rm (\ref{def:selector_for_beta_unknown})} is the best possible among all selectors
in model (\ref{model2}) in the \textit{asymptotically minimax} sense.
(As in publications \cite{BSt-2017} and \cite{ISt-2014}, the optimality of a selection
procedure here is understood in the minimax hypothesis testing sense.)

\begin{remark}\label{remark3}
	Inequalities {\rm (\ref{cond:inf})} and  {\rm (\ref{cond:inf2})}, which are analogous to inequalities {\rm(\ref{vector_cond:inf})} and {\rm (\ref{vector_cond:inf2})} in the vector case, describe the \textit{sharp selection boundary}, which defines a precise demarcation between what is possible and impossible in this problem.
	The sharp selection boundary is determined in terms of the function $a_{\varepsilon,k}(r_{\varepsilon,k})$ defined in {\rm (\ref{def:a2})},
	with $\mathring{\mathbb{Z}}^k$ instead of $\mathring{\mathbb{Z}}_u$,
	whose sharp asymptotics for every fixed $k$  are given by {\rm (\ref{aek})}.
	
\end{remark}

\begin{remark}\label{remark4}
	The {sharp selection boundary} determined by inequalities {\rm (\ref{cond:inf})} and  {\rm (\ref{cond:inf2})}
	is seen to lie above the sharp detection boundary as given  in {\rm (\ref{testing})}, with $a_{\varepsilon,k}(r_{\varepsilon,k})$ in place of
	$a_{\varepsilon,u}(r_{\varepsilon,k})$.
	This phenomenon, detected earlier for simpler models in \cite{BSt-2017} and \cite{ISt-2014}, is not unusual because, statistically,
	the problem of selecting nonzero components is more difficult than that of testing whether such nonzero components exist.
	In particular, in order to be \textit{selectable}, the components need to be \textit{detectable}.
\end{remark}

\subsection{The case of growing $k$}
We now state the analogues of Theorems \ref{th:UB_fixedK} and \ref{theorem:LB_figed_k} for the case $k=k_\varepsilon\to \infty$.
For this, the rates at which $d$ and $k$ tend to infinity need to be carefully regulated.

The upper bound on the maximum Hamming risk of our selector $\hat{\etab}$ is as follows.

\begin{theorem} \label{th:UB_growingK}
	Let $\beta\in(0,1)$ and $\sigma>0$ be fixed numbers, and let $d=d_\varepsilon \to \infty$ and $k=k_\varepsilon \to \infty$ be such that
	$k = \smallO(d)$, $\log \log d = \smallO (k)$, and $\log {d\choose k}=\smallO \left( \log \varepsilon^{-1} \right)$ as $\varepsilon\to 0$.
	Let the quantity $r_{\varepsilon,k}>0$ be as in Theorem {\rm \ref{th:UB_fixedK}}.
	Then the selector $\boldsymbol{\hat{\eta}}=(\hat{\eta}_u)_{u \in {\cal U}_{k,d} }$ given by {\rm (\ref{def:selector_for_beta_unknown})} 
	satisfies
	\begin{equation*}
		\limsup_{\varepsilon\to 0}\sup_{\boldsymbol{\eta} \in \mathcal{H}^d_{k,\beta}} \sup_{\thetab\in \Theta_{k,\sigma}^d(r_{\varepsilon,k})} \E{\thetab,\etab} |\boldsymbol{\eta} - \boldsymbol{\hat{\eta}}| = 0.
	\end{equation*}
\end{theorem}

The next result provides a lower bound on the minimax Hamming risk.
It shows that, under the same restrictions on  $d$ and $k$ as in Theorem \ref{th:UB_growingK}, if $r_{\varepsilon,k}$ is too ``small'' and
the identifiability condition
(\ref{cond:inf}) is not met, exact variable selection is impossible, and thus the selector
$\boldsymbol{\hat{\eta}}=(\hat{\eta}_u)_{u \in {\cal U}_{k,d} }$ given by {\rm (\ref{def:selector_for_beta_unknown})} is the best possible
in the asymptotically minimax sense with respect to the Hamming risk.

\begin{theorem} \label{theorem:LW_growing_k}
	Let $\beta\in(0,1)$ and $\sigma>0$ be fixed numbers, and let $d=d_\varepsilon \to \infty$ and $k=k_\varepsilon \to \infty$ be such that 
	$k = \smallO(d)$, $\log \log d = \smallO (k)$, and $\log {d\choose k}=\smallO \left( \log \varepsilon^{-1} \right)$ as $\varepsilon\to 0$.
	Let the quantity $r_{\varepsilon,k}>0$ be as in Theorem {\rm \ref{theorem:LB_figed_k}}.
	Then
	\begin{equation*}
		\liminf_{\varepsilon \to 0}\inf_{\tilde{\etab}} \sup_{\boldsymbol{\eta} \in \mathcal{H}^d_{k,\beta}} \sup_{\thetab\in \Theta_{k,\sigma}^d(r_{\varepsilon,k})} \E{\thetab,\etab} |\boldsymbol{\eta} - \boldsymbol{\tilde{\eta}}|>0,
	\end{equation*}
	where the infimum is taken over all selectors $\tilde{\etab} = (\tilde{\eta}_u)_{u \in \mathcal{U}_{k,d}}$ of a vector $\etab$ in model {\rm (\ref{model2})}. 
\end{theorem}

\begin{remark} \label{remark4}
	As seem from Theorems {\rm \ref{th:UB_growingK}} and {\rm \ref{theorem:LW_growing_k}}, for $k=k_\varepsilon\to \infty$,
	the sharp selection boundary given by  {\rm (\ref{cond:inf})} and {\rm (\ref{cond:inf2})}
	formally remains the same as for fixed $k$.
	It is determined in terms of the function $a_{\varepsilon,k}(r_{\varepsilon,k})$ defined in {\rm (\ref{def:a2})},
	with $\mathring{\mathbb{Z}}^k$ instead of $\mathring{\mathbb{Z}}_u$,
	whose sharp asymptotics, when  $k=k_\varepsilon\to \infty$, are given by relation {\rm (\ref{aek1})}.
\end{remark}

\section{Simulation study}
In this section, we illustrate numerically how the exact selector $\hat{\etab}=(\hat{\eta}_u)_{u\in {\cal U}_{k,d}}$ in (\ref{def:selector_for_beta_unknown})
works for identifying the active components of a signal $f\in {\cal F}_{k,\sigma}^d(r_{\varepsilon,k})$
of the form $f(\tb)=\sum_{u\in {\cal U}_{k,d}}\eta_u f_u(\tb_u)$ for $k=2$, $\sigma=1$, and $d=10, 50$.
The sparsity parameter $\beta$ and the cardinality ${d \choose 2}$ of the set $\mathcal{U}_{2,d}$ for different values of $d$ are as follows:
\begin{center}
	\begin{tabular}{ l | c | c }
		\hline
		$d$ & $\beta$ & ${d \choose 2}$ \\			
		\hline
		10 & 0.395 & 45 \\
		50 & 0.676 & 1225 \\
		\hline
	\end{tabular}
\end{center}
The two values of $\beta$, which were computed by using (\ref{sparsitycond})
as $\beta = 1- \log\left({\sum\nolimits_{u \in \mathcal{U}_{2,d}} \eta_u}\right)/{\log {d \choose 2}}$,
cover the ``dense case'' ($0 < \beta \leq {1}/{2}$) and the ``sparse case'' (${1}/{2} < \beta < 1 $).

\begin{figure}
	\centering
	\begin{subfigure}[b]{0.32\textwidth}
		\includegraphics[width=\textwidth]{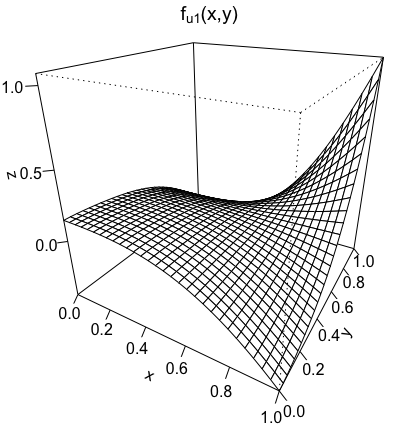}
	\end{subfigure}
	\begin{subfigure}[b]{0.32\textwidth}
		\includegraphics[width=\textwidth]{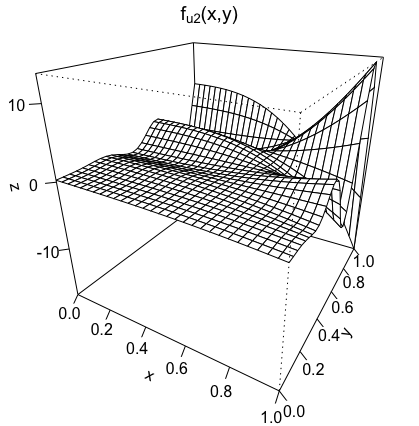}
	\end{subfigure}
	\begin{subfigure}[b]{0.32\textwidth}
		\includegraphics[width=\textwidth]{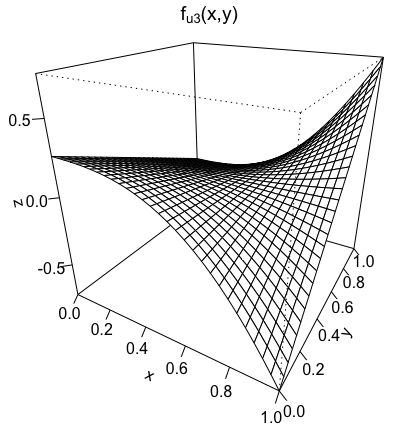}
	\end{subfigure}
	\begin{subfigure}[b]{0.32\textwidth}
		\includegraphics[width=\textwidth]{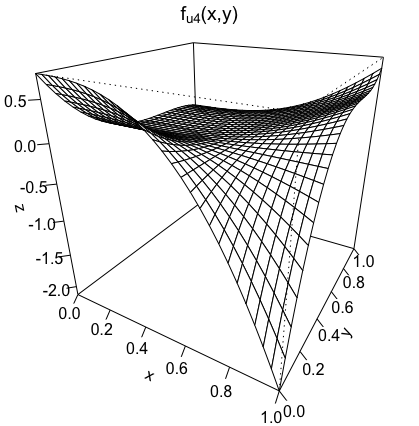}
	\end{subfigure}
	\begin{subfigure}[b]{0.32\textwidth}
		\includegraphics[width=\textwidth]{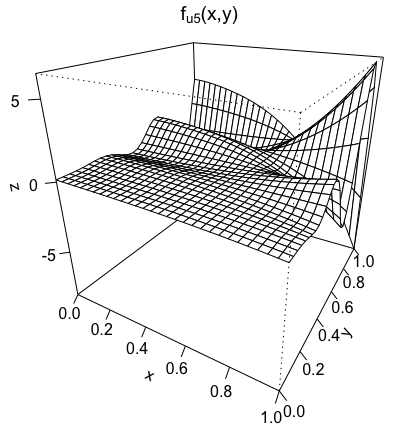}
	\end{subfigure}
	\begin{subfigure}[b]{0.32\textwidth}
		\includegraphics[width=\textwidth]{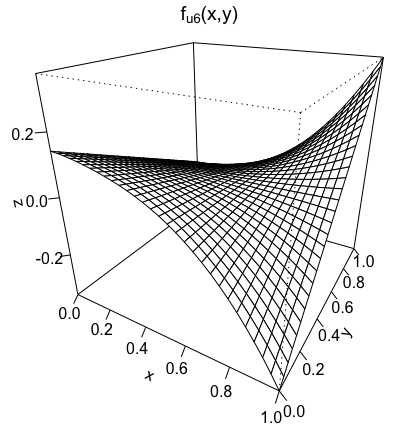}
	\end{subfigure}
	\begin{subfigure}[b]{0.32\textwidth}
		\includegraphics[width=\textwidth]{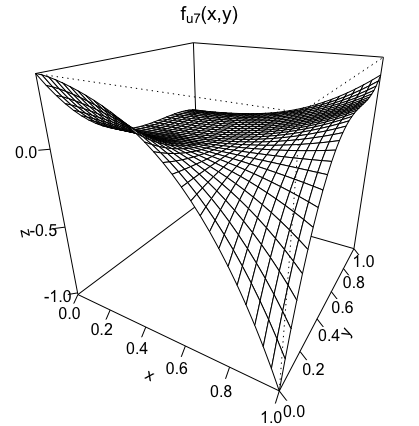}
	\end{subfigure}
	\begin{subfigure}[b]{0.32\textwidth}
		\includegraphics[width=\textwidth]{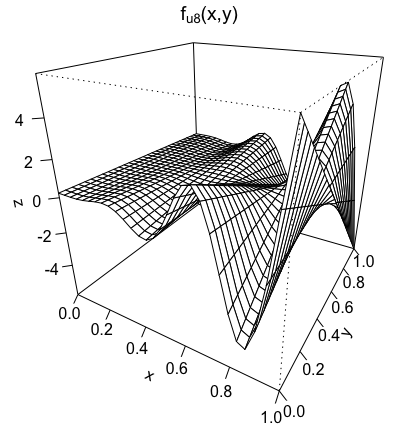}
	\end{subfigure}
	\begin{subfigure}[b]{0.32\textwidth}
		\includegraphics[width=\textwidth]{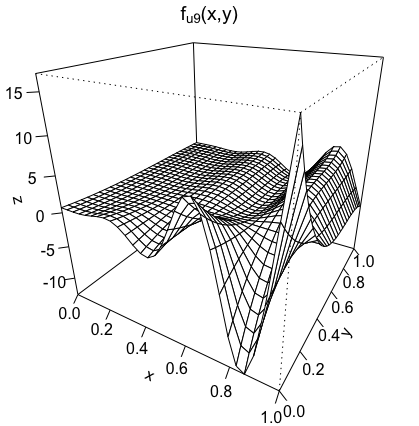}
	\end{subfigure}
	\begin{subfigure}[b]{0.32\textwidth}
		\includegraphics[width=\textwidth]{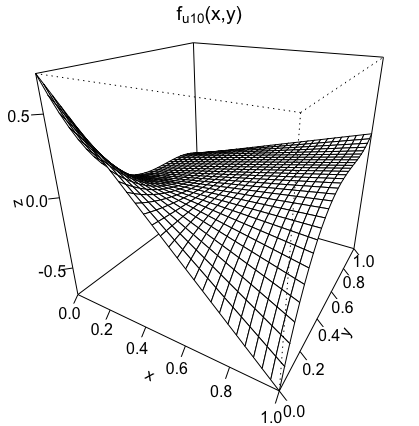}
	\end{subfigure}
	\caption{Bivariate components $f_u$ of the regression function $f$ defined by (\ref{def:functions_fu}).}\label{fig:fu}
\end{figure}

Consider the following five functions defined on $[0,1]$:
\begin{align*}
	g_1(t) &= t^2 (2^{t-1} - (t-0.5)^2 ) \exp(t) - 0.5424 , \\
	g_2(t) &= t^2 (2^{t-1} - (t-1)^5) - 0.2887, \\
	g_3 (t)&= 15 t^2 2^{t-1} \cos(15 t) - 0.5011, \\
	g_4 (t)&= t-{1}/{2}, \\
	g_5 (t)&= 5 (t-0.7)^3 + 0.29,
\end{align*}
and let the active components $f_u$ of $f$ defined on $[0,1]^2$ be as follows:
\begin{align}
	\begin{split}
		f_{u_1} (t_1,t_2)&= g_1(t_1) g_2(t_2),\qquad
		f_{u_2} (t_1,t_2)= g_1 (t_1)g_3(t_2), \\
		f_{u_3} (t_1,t_2)&= g_1 (t_1)g_4(t_2), \qquad
		f_{u_4} (t_1,t_2)= g_1 (t_1)g_5(t_2), \\
		f_{u_5} (t_1,t_2)&= g_2 (t_1)g_3(t_2),\qquad
		f_{u_6} (t_1,t_2)= g_2 (t_1)g_4(t_2), \\
		f_{u_7} (t_1,t_2)& = g_2 (t_1)g_5(t_2),\qquad
		f_{u_8} (t_1,t_2) = g_3 (t_1)g_4 (t_2), \\
		f_{u_9} (t_1,t_2) &= g_3 (t_1)g_5(t_2),\qquad
		f_{u_{10}} (t_1,t_2)= g_4 (t_1)g_5(t_2).
	\end{split}
	\label{def:functions_fu}
\end{align}
For every $u \in \mathcal{U}_{2,d}$,  the condition $\int_{0}^{1} f_u(\tb_u) dt_j = 0$ for all $j \in u,$ holds true up to four decimal places.
The bivariate functions defined in (\ref{def:functions_fu}) are depicted in Figure \ref{fig:fu}.
In our simulation study, we choose $u_1 = \{1,2\}$, $u_2 = \{1,3\}$, $u_3=\{1,4\}$, $u_4=\{1,5\}$, $u_5 =\{1,6\}$, $u_6=\{1,7\}$, $u_7=\{1,8\}$, $u_8=\{1,9\}$, $u_9=\{1,10\}$ for $d=10,50$.
We also take $u_{10}=\{2,1\}$ for $d=10$, and  $u_{10}=\{1,11\}$ for $d=50$. For~$u_i, \; i=11,\ldots,{d \choose 2}$, we set $f_{u_i} (t_1,t_2)=0.$

In order to simulate the observations in model (\ref{model2}), we need to compute the true Fourier coefficients  $\theta_\lb(u), \lb \in \mathring{\mathbb{Z}}_u$,
in decomposition (\ref{f-decomp}).
For $k=2$, we have $u=\{j_1, j_2\}$, $\tb_u = (t_{j_1},t_{j_2})$, and each index $\lb\in \mathring{\mathbb{Z}}_u$ has only two nonzero elements. We denote these nonzero elements of $\lb$ by $l_{j_1}$ and $l_{j_2}$. Then
\begin{align*}
	\theta_\lb(u) &= ( f_u, \phi_{l_{j_1}}  \phi_{l_{j_2}} )_{L_2^2}  = \int_{0}^{1} \int_{0}^{1} g_{j_1}(t_{j_1}) g_{j_2} (t_{j_2})  \phi_{l_{j_1}} (t_{j_1}) \phi_{l_{j_2}} (t_{j_2}) dt_{j_1} dt_{j_2} \nonumber \\
	&= \int_{0}^{1}  g_{j_1}(t_{j_1})   \phi_{l_{j_1}} (t_{j_1})  dt_{j_1} \int_{0}^{1}g_{j_2} (t_{j_2}) \phi_{l_{j_2}} (t_{j_2})dt_{j_2}  \nonumber \\
	& = ( g_{j_1} ,\phi_{l_{j_1}} )_{L^1_2}  (g_{j_2},  \phi_{l_{j_2}} )_{L^1_2}.
\end{align*}
For smooth Sobolev functions, the absolute values of their Fourier coefficients decay to zero at a~polynomial rate.
Therefore, although in theory $\lb=(l_1,\ldots,l_d) \in \mathring{\mathbb{Z}}_u$, we shall restrict ourselves to  $\lb \in \Ztu \deq  \mathring{\mathbb{Z}}_u \cap [-s,s]^d,$ where $s=3.6\times 10^4$ for $d=10$ and $s=3\times 10^4$ for $d=50.$
Thus, the model we are dealing with in this section  is as follows:
\begin{equation*}
	X_{\lb}=\eta_u\theta_\lb(u)+\varepsilon \xi_{\lb},\quad \lb \in \Ztu,
\end{equation*}
where $\varepsilon=10^{-4}$, the vector $(\xi_\lb)_{\lb \in \Ztu}$ consists of iid standard normal random variables $\xi_{\lb}$,  and the component $\eta_u $ of a vector $\etab=(\eta_u)_{u\in{\cal U}_{2,d}}\in{\cal H}_{2,\beta}^d$ equals 1 if $u \in \left\{ u_1,u_2,u_3, \ldots, u_{10}  \right\}$  and zero otherwise. 
Note that when $\varepsilon=10^{-4}$, $k=2,$ $\sigma=1,$ and $d=10,50,$ the condition $\log{d\choose k}<\log(\varepsilon^{-1})$ of Theorem~\ref{th:UB_fixedK} is satisfied.

Recall our definition of the exact selector $\hat{\etab}=(\hat{\eta}_u)_{u\in {\cal U}_{k,d}}$ in  (\ref{def:selector_for_beta_unknown}) with components 
\begin{equation*}
	\hat{\eta}_u = \max_{1 \leq m \leq M_k} \hat{\eta}_{u,m}, \quad  \hat{\eta}_{u,m} = \ind{ S_{u,m} > \sqrt{(2+\epsilon) \left( \log {d\choose k} + \log M_k \right)} },
\end{equation*}
where
$S_{u,m} = \sum_{\lb \in \mathring{\mathbb{Z}}_u} \omega_\lb (r^*_{\varepsilon,k,m}) \left[ \left( {X_\lb}/{\varepsilon}  \right)^2 -1 \right],$
$\omega_\lb (r^*_{\varepsilon,k,m}) = \dfrac{1}{2 \varepsilon^2} \dfrac{[\theta^*_\lb (r^*_{\varepsilon,k,m})]^2}{a_{\varepsilon,u}(r^*_{\varepsilon,k,m})},$
and the sharp asymptotics for $\theta^*_\lb (r_{\varepsilon,k})$ and $a_{\varepsilon,u}(r_{\varepsilon,k})$ are given by (\ref{thetal}) and (\ref{aek}).
We take $\epsilon=0.1$, $M_k=20$, and choose
the grid of points $\beta_{k,m}$, $m=1, \ldots, M_k$, uniformly on the interval $(0, 1)$.
For each $m$, the value $r^*_{\varepsilon,k,m}$ is found as a solution of equation (\ref{def:r_star_ekm}) with
$a_{\varepsilon,u} (r_{\varepsilon,k})$ as in (\ref{aek}).

For $d=10,50$, we run $J=12$ independent cycles of simulations and
estimate the Hamming risk ${\rm E}_{f,\etab}\left(\sum_{u\in{\cal U}_{2,d}}|\hat{\eta}_u-\eta_u|\right)$
by means of the quantity
$${Err}(\hat{\etab})=\frac{1}{J}\sum_{j=1}^J\sum_{u \in \mathcal{U}_{2,d}}|\hat{\eta}_u^{(j)}-\eta_u|,$$
where $\hat{\eta}_u^{(j)}$ is the value of $\hat{\eta}_u$ obtained in the $j$th repetition of the experiment.
The values of ${Err}(\hat{\etab})$ for different $d$ are seen in Table \ref{tab:Hast_err} in the column $\alpha=1$.

To examine the impact of the signal strength on the Hamming risk of $\hat{\etab}$, consider the function $f_{u_1}(t_1,t_2)=g_1(t_1)g_2(t_2)$
multiplied by $\alpha =0.01, 0.5, 1 ,2 ,5$.
The values of the
estimated Hamming risk ${Err}(\hat{\etab})$ obtained for different choices of $\alpha$ are presented in Table \ref{tab:Hast_err}.
\begin{table}[h!]
	\begin{center}
		\begin{tabular}{cc|c|c|c|c|c}
			&& \multicolumn{5}{ c }{$\alpha$} \\
			& & $0.01$ & $0.5$ & $1$ & $2$ & $5$ \\
			\hline
			{\multirow{2}{*}{$d$} } &  10 ($\beta = 0.40$)& 0.75 & 0.08 & 0 & 0  & 0   \\
			& 50 ($\beta = 0.68$) & 1.59 & 1.17 & 0.08 & 0  & 0 \\
		\end{tabular}
		\caption{Estimated Hamming risk of $\hat{\etab}$ obtained from 12 simulation cycles.}
		\label{tab:Hast_err}
	\end{center}
\end{table}

In general, as seen from Table \ref{tab:Hast_err}, the stronger the signal is, the smaller the number of misidentified active components of the signal $f$ is.
Note also that the sparser the model is, the harder the problem of selecting the active components of $f$ is.

\section{Proof of Theorems}

\textit{Proof of Theorem \ref{th:UB_fixedK}.} The proof goes partially along the lines of that of Theorem \ref{th:UB_fixedK} in \cite{ISt-2014}.
First, consider the $\chi^2$-type statistics $S_{k,m}$, $m=1,\ldots,M_k$, $1\leq k\leq d$, defined by
(\ref{def:S_um}) with $\mathring{\mathbb{Z}}^k$ instead of $\mathring{\mathbb{Z}}_u$,
and note some useful properties of these statistics.
The random variables
$(X_\lb/\varepsilon)^2$, $\lb\in \mathring{\mathbb{Z}}^k$, are independent and $(X_\lb/\varepsilon)^2\sim \chi_1^2((\theta_\lb/\varepsilon)^2)$
when
$\thetab_k=(\theta_\lb)_{\lb\in \mathring{\mathbb{Z}}^k} \in \mathring{\Theta}_{\cb_k}(r_{\varepsilon,k})$.
Therefore, noting that $\operatorname{E}(\chi_\nu^2(\lambda))=\nu+\lambda$, ${\rm var}(\chi_\nu^2(\lambda))=2(\nu+2\lambda)$,
and  taking into account (\ref{eq:sum_omega_is_half}), we obtain
\begin{align}
	\E{\thetab_k}(S_{k,m})&=\sum_{\lb\in\mathring{\mathbb{Z}}^k}\omega_\lb (r^*_{\varepsilon,k,m}) (\theta_\lb/\varepsilon)^2, \nonumber \\ %\label{Es}\\
	{\rm var}_{\thetab_k} (S_{k,m})& =\sum_{\lb \in \mathring{\mathbb{Z}}^k} {\rm var}_{\thetab_k} \left(\omega_\lb (r_{\varepsilon,k,m}^*) (X_\lb/\varepsilon)^2\right)=
	\sum_{\lb \in \mathring{\mathbb{Z}}^k} \omega^2_\lb (r_{\varepsilon,k,m}^*)\left(2+4(\theta_\lb/\varepsilon)^2\right)\nonumber\\
	&=
	1+4\sum_{\lb \in \mathring{\mathbb{Z}}^k} \omega^2_\lb (r_{\varepsilon,k,m}^*)(\theta_\lb/\varepsilon)^2
	\nonumber\\ &=  1+\mathcal{O}\left(\max_{\lb\in  \mathring{\mathbb{Z}}^k}\omega_\lb(r_{\varepsilon,k,m}^*)\E{\thetab_k}(S_{k,{m}})\right),\quad \varepsilon\to 0.\label{VarS}
\end{align}
It is clear that, under $\mathbf{H}_{0,k}$, we have
\begin{gather*}
	\operatorname{E}_{\boldsymbol{0}}(S_{k,m})=0,\quad {\rm{var}}_{\boldsymbol{0}}(S_{k,m})=1.
\end{gather*}
Also, let $T =T_{\varepsilon,k}\rightarrow \infty$ as $\varepsilon\to 0$ be such that
\begin{gather}
	T \max_{\lb \in \mathring{\mathbb{Z}}^k} \omega_\lb (r^*_{\varepsilon,k,m})=\smallO(1), \quad \varepsilon\to 0.\label{lem:ass_on_T}
\end{gather}
Then, for $m=1,\ldots,M_k,$ $1\leq k\leq d$, as $\varepsilon\to 0$
\begin{gather}\label{bound_1}
	\operatorname{P}_{\boldsymbol{0}} \left( S_{k,m} > T \right) \leq \exp \left( - \dfrac{T^2}{2}(1+\smallO(1)) \right).
\end{gather}
If, in addition to (\ref{lem:ass_on_T}), one has for $\thetab_k\in \mathring{\Theta}_{\cb_k}(r_{\varepsilon,k})$  and $m=1,\ldots,M_k,$ $1\leq k\leq d$,
\begin{gather}
	\E{\thetab_k} (S_{k,m}) \max_{\lb \in \mathring{\mathbb{Z}}^k} \omega_\lb (r^*_{\varepsilon,k,m}) =\smallO(1),
	\label{lem:ass_on_ES}
\end{gather}
then for this $\thetab_k$ and $m=1,\ldots,M_k,$ $1\leq k\leq d$, as $\varepsilon\to 0$
\begin{gather}\label{bound_2}
	\operatorname{P}_{\thetab_k} \left( S_{k,m} - \operatorname{E}_{\thetab_k}
	(S_{k,m}) \leq -T \right) \leq \exp \left( -\frac{T^2}{2}(1+\smallO(1)) \right)
\end{gather}
The proof of (\ref{bound_1}) and (\ref{bound_2}) is similar to that of Proposition 7.1 in \cite{DSA-2012} and therefore is omitted.

Next, for $r_{\varepsilon,k}>0$ consider the weights (see {\rm (\ref{def:omega})})
$$\omega_\lb (r_{\varepsilon,k})=\frac{1}{2\varepsilon^2}\frac{[\theta^*_{\lb}(r_{\varepsilon,k})]^2    }{a_{\varepsilon,k}(r_{\varepsilon,k})}, \quad \lb\in \mathring{\mathbb{Z}}^k.$$
Relations (\ref{thetal}) and (\ref{aek}) imply that for all those $\lb \in \mathring{\mathbb{Z}}^k$ for which the weights
$\omega_{\lb}(r_{\varepsilon,k})$ do not vanish,
\begin{gather}\label{omegal}
	\omega_{\lb}(r_{\varepsilon,k})\asymp r_{\varepsilon,k}^{k/(2\sigma)},\quad \varepsilon\to 0.
\end{gather}
The hypotheses $\mathbf{H}_{0,k}$ and $\mathbf{H}^{\varepsilon}_{1,k}$  separate asymptotically when
$r_{\varepsilon,k}/r^*_{\varepsilon,k}\to \infty$,
where, on one hand (see (\ref{def:r_star})),
$a_{\varepsilon,k}(r^*_{\varepsilon,k})\asymp \sqrt{\log {d\choose k}}$
and, on the other hand (see (\ref{aek})), $a_{\varepsilon,k}(r^*_{\varepsilon,k})\asymp \left(r^*_{\varepsilon,k}  \right)^{2+k/(2\sigma)}\varepsilon^{-2}$.
Therefore, using the ``continuity'' property of $a_{\varepsilon,k}(r_{\varepsilon,k})$,
as stated in (\ref{def:continuity_of_a}), for $m=1,\ldots,M_k$, $1\leq k\leq d$, we obtain, as $\varepsilon\to 0$
\begin{gather}\label{rek}
	r^*_{\varepsilon,k}\asymp r^*_{\varepsilon,k,m}\asymp\left(\varepsilon\left[\log {d\choose k}  \right]^{1/4}  \right)^{{4\sigma}/{(4\sigma+k)}   },
\end{gather}
and hence, by using (\ref{omegal}), the weights of the $\chi^2$-type statistics $S_{k,m}$, $m=1,\ldots,M_k$, $1\leq k\leq d$, satisfy
\begin{equation}\label{omegamax}
	\max_{\lb \in \mathring{\mathbb{Z}}^k} \omega_\lb(r_{\varepsilon,k,m}^*) \asymp \left( \varepsilon \left[ \log {d\choose k} \right]^{{1}/{4}} \right)^{{2k}/{(4 \sigma +k)}}.
\end{equation}

Now, consider the maximum risk of the selector $\boldsymbol{\hat{\eta}}=(\hat{\eta}_u)_{u\in {\cal U}_{k,d}}$. We have
\begin{align}
	R_{\varepsilon,k}(\boldsymbol{\hat{\eta}}) :=& \sup_{\boldsymbol{\eta} \in \mathcal{H}^d_{k,\beta}} \sup_{\thetab \in \Theta_{k,\sigma}^d(r_{\varepsilon,k}) } \E{\thetab,\etab} |\boldsymbol{\eta} - \boldsymbol{\hat{\eta}}| = \sup_{\boldsymbol{\eta} \in \mathcal{H}^d_{k,\beta}} \sup_{\thetab \in \Theta_{k,\sigma}^d(r_{\varepsilon,k})} \E{\thetab,\etab} \sum_{u\in{\cal U}_{k,d}} |\eta_u - \hat{\eta}_u| \nonumber \\
	=& \sup_{\boldsymbol{\eta} \in \mathcal{H}^d_{k,\beta}} \sup_{\thetab \in \Theta_{k,\sigma}^d(r_{\varepsilon,k})} \E{\thetab,\etab} \left( \sum_{\substack{u\,:\, \eta_u=0}} \hat{\eta}_u + \sum_{\substack{u\,:\, \eta_u=1}} (1-\hat{\eta}_u)  \right). \label{eq:R_part2}
\end{align}
Using the notation 
\begin{equation*}
t(\varepsilon,k)=\sqrt{(2+\epsilon) \left( \log {d\choose k} + \log M_k \right)},  
\end{equation*}
this can be estimated from above by
\begin{align*}
	R_{\varepsilon,k}(\boldsymbol{\hat{\eta}})  & \leq \sup_{\boldsymbol{\eta} \in \mathcal{H}^d_{k,\beta}} \left( \sum_{\substack{u\,: \, \eta_u=0}} \E{{\bf 0},\eta_u} (\hat{\eta}_u) + \sum_{\substack{u\,: \, \eta_u=1}} \sup_{\thetab_u \in \mathring{\Theta}_{\cb_u}(r_{\varepsilon,k})} \E{\thetab_u,\eta_u} (1-\hat{\eta}_u)  \right) \nonumber \\
	&= \sup_{\boldsymbol{\eta} \in \mathcal{H}^d_{k,\beta}} \left( \sum_{\substack{u\,:\, \eta_u=0}} \E{{\bf 0},\eta_u} \left\{ \max_{1 \leq m \leq M_k} \ind{ S_{u,m} >  t(\varepsilon,k)} \right\} + \right.\nonumber \\
	& \quad \left. +  \sum_{\substack{u\,: \, \eta_u=1}} \sup_{\thetab_u \in \mathring{\Theta}_{\cb_u}(r_{\varepsilon,k})} \E{\thetab_u,\eta_u}
	\left\{ 1-\max_{1 \leq m \leq M_k} \ind{ S_{u,m} > t(\varepsilon,k) }  \right\}  \right).
\end{align*}
From this, noting that $\sum_{\substack{u\,: \, \eta_u=1}} \eta_u = {d\choose k}^{1-\beta}(1 + \smallO(1)) \leq 2 {d\choose k}^{1-\beta}$
and switching from $S_{u,m}$ to $S_{k,m}$, we obtain for all small enough $\varepsilon$
\begin{align}
	R_{\varepsilon,k}(\boldsymbol{\hat{\eta}}) & \leq {d\choose k}  \E{{\bf 0}} \left[ \max_{1 \leq m \leq M_k} \ind{ S_{k,m} > t(\varepsilon,k) } \right] + \nonumber \\
	& \quad + 2 {d \choose k}^{1-\beta} \sup_{\thetab_k \in \mathring{\Theta}_{\cb_k}(r_{\varepsilon,k})} \E{\thetab_k} \left\{ 1-\max_{1 \leq m \leq M_k} \ind{ S_{k,m} > t(\varepsilon,k)}  \right\} \nonumber\\
	& =: J_{\varepsilon,k}^{(1)} + J_{\varepsilon,k}^{(2)},
	\label{risk}
\end{align}

For the  term $J_{\varepsilon,k}^{(1)}$,
in view of (\ref{omegamax}) and the assumptions $\log {d \choose k}=\smallO \left( \log \varepsilon^{-1} \right)$ and
$\log M_k=\smallO(\log {d\choose k})$,  as $\varepsilon\to 0$
\begin{gather*}
	t(\varepsilon,k)	\max_{\lb \in \mathring{\mathbb{Z}}^k} \omega_\lb(r^*_{\varepsilon,k,m_0})  
	\asymp \sqrt{\log {d\choose k}}\left( \varepsilon \left[ \log {d\choose k} \right]^{{1}/{4}} \right)^{{2k}/{(4 \sigma +k)}}
	= \smallO(1),
\end{gather*}
and thus (\ref{lem:ass_on_T}) with $T=t(\varepsilon,k)$ is satisfied. Therefore, by using (\ref{bound_1}) with $T=t(\varepsilon,k)$,
as $\varepsilon\to 0$
\begin{align}
	J_{\varepsilon,k}^{(1)}	 &\leq {d\choose k} \sum_{m=1}^{M_k} \E{{\bf 0}} \left\{ \ind{ S_{k,m} > t(\varepsilon,k) } \right\} \leq M_k {d\choose k} \max_{1\leq m\leq M_k}\Prob{{\bf 0}}{S_{k,m} > t(\varepsilon,k)}\nonumber\\
	&\leq M_k {d\choose k} \exp \left( -\frac{2+\epsilon}{2} \left( \log {d\choose k} +\log M_k\right)  \right) 
	\asymp	\left( M_k {d\choose k}\right)^{-\frac{\epsilon}{2} } = \smallO(1),
	\label{eq:ThUB_I1_final}
\end{align}
where the last equality holds since $\log M_k=\smallO(\log {d\choose k})$ and condition (\ref{delta}).
For the term $J_{\varepsilon,k}^{(2)}$ in (\ref{risk}), we have
\begin{align}
	J_{\varepsilon,k}^{(2)}	 &= 2{d \choose k}^{1-\beta} \sup_{\thetab_k \in \mathring{\Theta}_{\cb_k}(r_{\varepsilon,k})} \E{\thetab_k} \left[ \ind{ \bigcap\limits_{1 \leq m \leq M_k} \left\{ S_{k,m} \leq t(\varepsilon,k) \right\} }  \right] \nonumber \\
	&=2{d \choose k}^{1-\beta} \sup_{\thetab_k \in \mathring{\Theta}_{\cb_k}(r_{\varepsilon,k})} \operatorname{P}_{\thetab_k} \left( \bigcap\limits_{1 \leq m \leq M_k} \left\{ S_{k,m} \leq t(\varepsilon,k) \right\}   \right) \nonumber \\
	&\leq2 {d \choose k}^{1-\beta} \sup_{\thetab_k \in \mathring{\Theta}_{\cb_k}(r_{\varepsilon,k})} \min_{1 \leq m \leq M_k} \operatorname{P}_{\thetab_k} \left(   S_{k,m} \leq t(\varepsilon,k)\right).
	\label{eq:ThUB_I2_2}
\end{align}

To obtain a good upper bound on the right-hand side of (\ref{eq:ThUB_I2_2}), we shall use the result that is similar to Proposition 4.1 in \cite{DSA-2012}:
for $B \geq 1, \varepsilon>0$, and $r_{\varepsilon,k}>0$, $1\leq k\leq d$,
\begin{equation}\label{lemma:ball}
	\frac{1}{\varepsilon^2} \inf_{\thetab_k \in \mathring{\Theta}_{\cb_k}(Br_{\varepsilon,k})} \sum_{\lb \in \mathring{\mathbb{Z}}^k} \omega_\lb (r_{\varepsilon,k}) \theta_\lb^2 \geq B^2 a_{\varepsilon,k}(r_{\varepsilon,k}).
\end{equation}
First, recall the equidistant grid points $\beta_{k,m}=m\rho_k$, $m=1,\ldots,M_k$, on the interval $(0,1)$ and let index $m_0$ be such that  $\beta\in(\beta_{k,m_0},\beta_{k,m_0+1}]$.
Condition  (\ref{cond:inf}) implies that for all small enough $\varepsilon$ and some $\Delta_0>0$ ($\Delta_0\to 0 $ as $\varepsilon\to 0$),
\begin{gather*}
	a_{\varepsilon,k}(r_{\varepsilon,k})\geq \sqrt{2}(1+\sqrt{1-\beta})\sqrt{\log {d\choose k}}(1+\Delta_0).
\end{gather*}
Next, by (\ref{def:r_star_ekm}) and the choice of $m_0$, for all small enough $\varepsilon$ and some $\Delta_1>0$ ($\Delta_1$ may be taken to have $\Delta_1<\Delta_0/3$),
\begin{gather*}
	a_{\varepsilon,k}(r^*_{\varepsilon,k,m_0})=\sqrt{2}(1+\sqrt{1-\beta})\sqrt{\log {d\choose k}}(1+\Delta_1).
\end{gather*}
From this, by the ``continuity'' property of $a_{\varepsilon,k}$ as presented in (\ref{def:continuity_of_a}),
for all small enough $\varepsilon$ and some (small) $\Delta_2, \Delta_3>0$ (we may take $\Delta_3<\Delta_0/3$), we obtain
\begin{align*}
	a_{\varepsilon,k}((1+\Delta_2)r^*_{\varepsilon,k,m_0})&\leq (1+\Delta_3)a_{\varepsilon,k}(r^*_{\varepsilon,k,m_0}) = \sqrt{2}(1+\sqrt{1-\beta})\sqrt{\log {d\choose k}}(1+\Delta_1)(1+\Delta_3)\\
	&\leq \sqrt{2}(1+\sqrt{1-\beta})\sqrt{\log {d\choose k}}(1+\Delta_0)\leq a_{\varepsilon,k}(r_{\varepsilon,k}),
\end{align*}
and hence, by the monotonicity of $a_{\varepsilon,k}(r_{\varepsilon,k})$,
\begin{gather*}
	r_{\varepsilon,k}\geq (1+\Delta_2)r^*_{\varepsilon,k,m_0}=:Br^*_{\varepsilon,k,m_0},\quad B>1.
\end{gather*}
Using (\ref{lemma:ball}), we get
for some $\Delta_4>0$ such that $B^2 (1+\Delta_1) > 1+\Delta_4$  for all small enough $\varepsilon$ that
\begin{align}
	\inf_{\thetab_k \in \mathring{\Theta}_{\cb_k}(r_{\varepsilon,k})} \E{\thetab_k}(S_{k,m_0}) &\geq\varepsilon^{-2}
	\inf_{\thetab_k \in \mathring{\Theta}_{\cb_k}(Br_{\varepsilon,k,m_0}^*)} \sum_{\lb \in \mathring{\mathbb{Z}}^k} \omega_\lb(r_{\varepsilon,k,m_0}^*) \theta_\lb^2 \geq B^2 a_{\varepsilon,k} (r_{\varepsilon,k,m_0}^*)  \nonumber \\
	&= B^2 \sqrt{2}(1+\sqrt{1-\beta})\sqrt{\log {d\choose k}}(1+\Delta_1) \nonumber \\
	&> \sqrt{2}(1+\sqrt{1-\beta})\sqrt{\log {d\choose k}}(1+\Delta_4).
	\label{cor_bound_infES}
\end{align}
We may choose $\Delta_4=\smallO(1)$ to be such that $\Delta_4\log {d\choose k}\to \infty$ and $\epsilon=\smallO(\Delta_4)$ as $\varepsilon\to 0.$
Since by assumption $\log M_k = \smallO \left( \log {d\choose k} \right)$ as $\varepsilon \to 0$, it follows from (\ref{cor_bound_infES}) that
for all small enough $\varepsilon$ 
\begin{multline*}
	t(\varepsilon,k) - \inf_{\thetab_k \in \mathring{\Theta}_{\cb,k}(r_{\varepsilon,k})} \E{\thetab_k}( S_{k,{m_0}}) \leq  \left( \sqrt{2+\epsilon} - \sqrt{2}(1+\sqrt{1-\beta})(1+\Delta_4) \right) \sqrt{\log {d\choose k}}(1+\smallO(1)) =:-\mathbb{T},
\end{multline*}
where $\mathbb{T}=\mathbb{T}_{\varepsilon,k}\to \infty$ as $\varepsilon\to 0.$

Using this fact, we shall now estimate the term on the right-hand side of (\ref{eq:ThUB_I2_2}). 
For this, letting the  index $m_0$ be chosen as before, define the subsets $\mathring{\Theta}^{(i)}_{\cb_k,m_0}(r_{\varepsilon,k})$, $i=1,2,3$, of
$\mathring{\Theta}_{\cb_k}(r_{\varepsilon,k})$ as follows:
\begin{eqnarray*}
	\mathring{\Theta}^{(1)}_{\cb_k,m_0}(r_{\varepsilon,k})\!\!\!&=&\!\!\!\left\{\thetab_k\in \mathring{\Theta}_{\cb_k}(r_{\varepsilon,k}):
	\limsup_{\varepsilon\to 0} \E{\thetab_k}(S_{k,m_0}) \max_{\lb \in \mathring{\mathbb{Z}}^k} \omega_\lb (r^*_{\varepsilon,k,m_0})=0 \right\},\\
	\mathring{\Theta}^{(2)}_{\cb_k,m_0}(r_{\varepsilon,k})\!\!\!&=&\!\!\!\left\{\thetab_k\in \mathring{\Theta}_{\cb_k}(r_{\varepsilon,k}):c\leq \liminf_{\varepsilon\to 0} \E{\thetab_k}(S_{k,m_0}) \max_{\lb \in \mathring{\mathbb{Z}}^k} \omega_\lb (r^*_{\varepsilon,k,m_0})\leq\right. \\
	& &\quad\quad\leq \limsup_{\varepsilon\to 0} \E{\thetab_k}(S_{k,m_0}) \max_{\lb \in \mathring{\mathbb{Z}}^k} \omega_\lb (r^*_{\varepsilon,k,m_0})  \leq C\; \mbox{for some}\; 0<c\leq C<\infty\bigg\},\\
	\mathring{\Theta}^{(3)}_{\cb_k,m_0}(r_{\varepsilon,k})\!\!\!&=&\!\!\!\left\{\thetab_k\in \mathring{\Theta}_{\cb_k}(r_{\varepsilon,k}):
	\liminf_{\varepsilon\to 0} \E{\thetab_k}(S_{k,m_0}) \max_{\lb \in \mathring{\mathbb{Z}}^k} \omega_\lb (r^*_{\varepsilon,k,m_0})=\infty \right\},
\end{eqnarray*}
and note that relation (\ref{omegamax}) holds true for all $1\leq m\leq M_k$, including $m_0$.
Since, $\mathring{\Theta}^{(i)}_{\cb_k,m_0}(r_{\varepsilon,k})$, $i=1,2,3$, form a partition of
$\mathring{\Theta}_{\cb_k}(r_{\varepsilon,k})$, we may continue from~(\ref{eq:ThUB_I2_2}):
\begin{align}
	& J_{\varepsilon,k}^{(2)}\leq
	2{d \choose k}^{1-\beta}\sum_{i=1}^3 \sup_{\thetab_k\in \mathring{\Theta}^{(i)}_{\cb_k,m_0}(r_{\varepsilon,k})}
	\operatorname{P}_{\thetab_k} \Bigg(   S_{k,m_0}
	\leq
	t(\varepsilon,k) 
	\Bigg)\nonumber \\
	&\leq 2{d \choose k}^{1-\beta}\sup_{\thetab_k\in \mathring{\Theta}^{(1)}_{\cb_k,m_0}(r_{\varepsilon,k})}
	\operatorname{P}_{\thetab_k} \Bigg(   S_{k,m_0}- \E{\thetab_k} (S_{k,m_0}) \leq -\mathbb{T}\Bigg)+ 2{d \choose k}^{1-\beta}\times\nonumber\\
	& \quad \times\sum_{i=2}^3\sup_{\thetab_k\in \mathring{\Theta}^{(i)}_{\cb_k,m_0}(r_{\varepsilon,k})}
	\operatorname{P}_{\thetab_k} \Bigg(   S_{k,m_0}- \E{\thetab_k} (S_{k,m_0}) \leq t(\varepsilon,k)  -\E{\thetab_k} (S_{k,m_0})\Bigg)
	\nonumber\\ 
	&=:
	K_{\varepsilon,k,m_0}^{(1)}+ K_{\varepsilon,k,m_0}^{(2)}+K_{\varepsilon,k,m_0}^{(3)}.
	\label{J2}
\end{align}

Consider the  term $ K_{\varepsilon,k,m_0}^{(1)}$.
If $\thetab_k \in \mathring{\Theta}^{(1)}_{\cb_k,m_0}(r_{\varepsilon,k})$, then condition (\ref{lem:ass_on_ES}) holds,
and the application of the upper bound (\ref{bound_2}) with $T=\mathbb{T}\asymp\sqrt{\log {d\choose k}  }$ as above gives
\begin{align}
	& K_{\varepsilon,k,m_0}^{(1)}
	\leq 2{d \choose k}^{1-\beta} \exp \left( -\frac{\mathbb{T}^2}{2} (1+\smallO(1))\right)  \nonumber \\
	&=2 {d \choose k}^{1-\beta} \exp \left( -\frac{( \sqrt{2+\epsilon} - \sqrt{2}(1+\sqrt{1-\beta})(1+\Delta_4))^2}{2} \log {d \choose k} (1+\smallO(1)) \right) \nonumber \\
	&\asymp {d \choose k}^{-\beta - \frac{\epsilon}{2} - (1+\sqrt{1-\beta})^2(1+\Delta_4)^2 + \sqrt{2(2+\epsilon)}
		(1+\sqrt{1-\beta})(1+\Delta_4) }=\smallO(1). 
	\label{K21}
\end{align}

When $\thetab_k$ falls in the sets $ \mathring{\Theta}^{(2)}_{\cb_k,m_0}(r_{\varepsilon,k})$ and
$ \mathring{\Theta}^{(3)}_{\cb_k,m_0}(r_{\varepsilon,k})  $, condition (\ref{lem:ass_on_ES})
does not hold anymore, and hence (\ref{bound_2}) is not applicable.
For deriving good upper bounds for the terms $K_{\varepsilon,k,m_0}^{(2)}$  and $K_{\varepsilon,k,m_0}^{(3)}$ in~(\ref{J2}), we shall apply Chebyshev's inequality.
First, by Chebyshev's inequality and relations  (\ref{VarS}) and~(\ref{omegamax}), 
\begin{align}
&	K_{\varepsilon,k,m_0}^{(2)}\leq 2 {d \choose k}^{1-\beta}  \sup_{\thetab_k\in \mathring{\Theta}^{(2)}_{\cb_k,m_0}(r_{\varepsilon,k})}
	\frac{{\rm var}_{\thetab_k} (S_{k,m_0})  }{\E{\thetab_k}^2(S_{k,m_0}) (1+\smallO(1) ) }\nonumber\\ &\leq
	2{d \choose k}^{1-\beta}  \sup_{\thetab_k\in \mathring{\Theta}^{(2)}_{\cb_k,m_0}(r_{\varepsilon,k})      }\frac{  1+\mathcal{O}\left(\max_{\lb\in  \mathring{\mathbb{Z}}^k}\omega_\lb(r_{\varepsilon,k,m_0}^*)\E{\thetab_k}(S_{k,{m_0}})\right)}{\E{\thetab_k}^2(S_{k,m_0})(1+\smallO(1))  }\nonumber\\
	&=\mathcal{O}\left({d \choose k}^{1-\beta}\left( \varepsilon \left[ \log {d\choose k} \right]^{{1}/{4}} \right)^{{4k}/{(4 \sigma +k)}} \right)=\smallO(1),
	\quad \varepsilon\to 0,
	\label{K22}
\end{align}
where the last equality holds since $\log{d\choose k}=\smallO(\log \varepsilon^{-1})$ as $\varepsilon\to 0$.
Next, for the term  $K_{\varepsilon,k,m_0}^{(3)}$, by Chebyshev's inequality and relations (\ref{VarS}), (\ref{omegamax}), and (\ref{cor_bound_infES}), we have as $\varepsilon\to 0$
\begin{align}
	K_{\varepsilon,k,m_0}^{(3)}&\leq 2 {d \choose k}^{1-\beta}  \sup_{\thetab_k\in \mathring{\Theta}^{(3)}_{\cb_k,m_0}(r_{\varepsilon,k})}
	\frac{{\rm var}_{\thetab_k} (S_{k,m_0})  }{\E{\thetab_k}^2(S_{k,m_0}) (1+\smallO(1) ) } \nonumber \\
	&\leq\frac{{d \choose k}^{1-\beta}\mathcal{O}\left(\max_{\lb\in  \mathring{\mathbb{Z}}^k}\omega_\lb(r_{\varepsilon,k,m_0}^*)\right)}{\inf_{\thetab_k\in \mathring{\Theta}^{(3)}_{\cb_k,m_0}(r_{\varepsilon,k})}\E{\thetab_k}(S_{k,m_0})}\nonumber \\
	&= \mathcal{O}\left({d \choose k}^{1-\beta} \left[\log{d\choose k}\right]^{-1/2}\left( \varepsilon \left[ \log {d\choose k} \right]^{{1}/{4}} \right)^{{2k}/{(4 \sigma +k)}} \right) =\smallO(1),\label{K23}
\end{align}
where the last equality holds since $\log{d\choose k}=\smallO(\log \varepsilon^{-1})$ as $\varepsilon\to 0$.

Combining (\ref{risk}), (\ref{eq:ThUB_I1_final}) and (\ref{J2}) to (\ref{K23}),
we obtain for all $\beta\in(0,1)$ and $\sigma>0$
$$R_{\varepsilon,k}(\boldsymbol{\hat{\eta}}) \leq J^{(1)}_{\varepsilon,k} + J^{(2)}_{\varepsilon,k} = \smallO(1),\quad \varepsilon\to \infty,$$
which completes the proof.  \qedsymbol
\medskip

\textit{Proof of Theorem \ref{theorem:LB_figed_k}.}
Since
$$\sup_{\thetab\in \Theta_{k,\sigma}^d(r^{\prime}_{\varepsilon,k})} \E{\thetab,\etab} |\boldsymbol{\eta} - \boldsymbol{\tilde{\eta}}|\geq
\sup_{\thetab\in \Theta_{k,\sigma}^d(r^{\prime\prime}_{\varepsilon,k})} \E{\thetab,\etab} |\boldsymbol{\eta} - \boldsymbol{\tilde{\eta}}|$$
whenever $0<r^{\prime}_{\varepsilon,k}<r^{\prime\prime}_{\varepsilon,k}$,
by the monotonicity of $a_{\varepsilon,k}$, we can restrict ourselves to the case when
\begin{gather}\label{adas}
	\liminf_{\varepsilon \to 0} \frac{a_{\varepsilon,k}(r_{\varepsilon,k})}{\sqrt{\log {d\choose k}}} >0,%\leq	
\end{gather}
which, together with (\ref{cond:inf2}), gives $a_{\varepsilon,k}(r_{\varepsilon,k})\asymp \sqrt{\log {d\choose k}}$ as $\varepsilon\to 0.$

For  $1\leq k\leq d$, suppressing for brevity the dependence on $k$, denote by $p = {d \choose k}^{-\beta}$
the proportion of nonzero components of $\etab = (\eta_u)_{u \in \mathcal{U}_{k,d}},$ and introduce the prior distribution of $\etab$ (cf. the prior distribution in Section 7.3 of \cite{DSA-2012}) as follows:
\begin{equation*}
	\pi_{\etab} = \prod_{u \in \mathcal{U}_{k,d}} \pi_{\eta_u}, \quad  \pi_{\eta_u}=(1-p) \delta_0 + p \delta_1,
	\label{LB:prior_distr_eta}
\end{equation*}
where $\delta_x$ is the $\delta$-measure that puts a pointmass 1 at $x$.
That is,  $\eta_{u}$, $u\in {\cal U}_{k,d}$, are independent Bernoulli random variables with parameter $p$. Further,
recall the sequence $(\theta^*_\lb (r_{\varepsilon,k}))_{\lb \in \mathring{\mathbb{Z}}_u}$ as in~(\ref{def:a2_with_theta_star}), let
$\theta_\lb^*=\theta^*_{\lb}(r_{k,\varepsilon})$,
and define the prior distribution of $\thetab = (\thetab_{u})_{u \in \mathcal{U}_{k,d}}$ by
\begin{equation*}
	\pi_{\thetab} = \prod_{u \in \mathcal{U}_{k,d}} \pi_{\thetab_{u}}, \quad  \pi_{\thetab_{u}} = \prod_{\lb \in \mathring{\mathbb{Z}}_u} \frac{\delta_{\theta_\lb^*} + \delta_{-\theta_\lb^*}}{2}.
	\label{LB:prior_distr_thetab}
\end{equation*}
The random variables $\theta_{\lb}$, $\lb\in\mathring{\mathbb{Z}}_u$, $u\in {\cal U}_{k,d}$, are independent, and each $\theta_\lb$ takes on the values $\theta_\lb^*$ and $-\theta_\lb^*$ with probability 1/2. By setting
$$R_{\varepsilon,k} := \inf_{\tilde{\etab}} \sup_{\boldsymbol{\eta} \in \mathcal{H}^d_{k,\beta}} \sup_{\thetab\in \Theta_{k,\sigma}^d(r_{\varepsilon,k})} \E{\thetab,\etab} |\boldsymbol{\eta} - \boldsymbol{\tilde{\eta}}|,$$
we can estimate the minimax risk $R_{\varepsilon,k}$ from below by the Bayes risk as follows:
\begin{align}
	R_{\varepsilon,k} &\geq \inf_{\tilde{\etab}} \E{\pi_{\etab}} \E{\pi_{\thetab}} \operatorname{E}_{\etab,\thetab} |\etab - \tilde{\etab}| = \inf_{\tilde{\etab}} \E{\pi_{\etab}} \E{\pi_{\thetab}} \operatorname{E}_{\etab,\thetab} \left( \sum_{u \in \mathcal{U}_{k,d}} |\eta_{u} - \tilde{\eta}_u| \right) \nonumber \\
	&= \inf_{\tilde{\etab}} \sum_{u \in \mathcal{U}_{k,d}} \E{\pi_{\eta_u}} \E{\pi_{\thetab_{u}}} \operatorname{E}_{\eta_u,\thetab_{u}} |\eta_{u} - \tilde{\eta}_u|.
	\label{LB:risk_uneq1}
\end{align}
Consider the data $\Xb_{\!u} = (X_\ell)_{\ell \in \mathring{\mathbb{Z}}_u}$, $u\in {\cal U}_{k,d}$, where  $X_\lb \sim N (\eta_u \theta_\lb,\varepsilon^2)$, and
introduce the following continuous mixture of distributions:
\begin{equation*}
	\operatorname{P}_{\pi,\eta_u} (d\Xb_{\!u}) = \E{\pi_{\thetab_{u}}} \operatorname{P}_{\eta_{u},\thetab_{u}} (dX_\ell),\quad u\in {\cal U}_{k,d},
	\label{LB:mixture_distr_of_thetab_and_X}
\end{equation*}
that is,
\begin{equation*}
	\E{\pi_{\thetab_{u}}} \operatorname{E}_{\eta_u,\thetab_{u}} |\eta_{u} - \tilde{\eta}_u| = \E{\pi_{\thetab_{u}}} \int |\eta_{u} - \tilde{\eta}_u| d\operatorname{P}_{\eta_u,\thetab_{u}} = \int |\eta_{u} - \tilde{\eta}_u| d\operatorname{P}_{\pi,\eta_u}.
\end{equation*}
This mixture of distributions can be alternatively expressed as
\begin{equation*}
	\operatorname{P}_{\pi,\eta_{u}} = \prod_{\lb \in \mathring{\mathbb{Z}}_u} \left[ \frac{N(\eta_u \theta_\lb^*,\varepsilon^2) + N(- \eta_u \theta_\lb^*,\varepsilon^2)}{2} \right], \quad u \in \mathcal{U}_{k,d}.
\end{equation*}
Set
\begin{equation*}
	\nu_{\lb}^* := {\theta_\lb^*}/{\varepsilon},
	\label{def:nu_star2}
\end{equation*}
and define the independent random variables
\begin{equation*}
	Y_\lb := \frac{X_\lb}{\varepsilon} = \eta_{u} \nu_\lb^* +  \xi_{\lb}
	\sim N(\eta_{u} \nu_\lb^*,1), \quad \lb \in \mathring{\mathbb{Z}}_u,\quad u \in \mathcal{U}_{k,d},
	\label{LB:Yl_normal}
\end{equation*}
where the dependence of  $\nu_{\lb}^*$ on $u$ is omitted.
Then, denoting ${\Yb}_{\!\!u} = (Y_\lb)_{\lb \in \mathring{\mathbb{Z}}_u}$, we can express the likelihood ratio in the form
\begin{equation}
	\Lambda_{\pi,u} := \frac{d \operatorname{P}_{\pi,1}}{d \operatorname{P}_{\pi,0}} (\Yb_{\!\!u}) = \prod_{\lb \in \mathring{\mathbb{Z}}_u} \exp \left( - \frac{(\nu_{\lb}^*)^2}{2} \right) \cosh \left( (\nu_{\lb}^*)^2 Y_\lb \right),
	\label{LB:likelihood_ratio}
\end{equation}
where $\nu_{\lb}^*$, $\lb \in \mathring{\mathbb{Z}}_u$, satisfy as $\varepsilon\to 0$
\begin{gather}\label{LB:nu_star_is_o1}
	\nu_{\lb}^*=\smallO(1).
\end{gather}
Indeed, recall the quantity $r^*_{\varepsilon,k}>0$ determined by (\ref{def:r_star}).
By  (\ref{cond:inf2}) and the ``continuity'' of $a_{\varepsilon,u}$, it holds that
$r_{\varepsilon,k}/r^*_{\varepsilon,k} <1$ for all small enough $\varepsilon$.
Next, by (\ref{thetal}) and (\ref{rek}), as $\varepsilon\to 0$
\begin{gather*}
	(\nu_{\lb}^*)^2 \asymp \varepsilon^{-2} r_{\varepsilon,k}^{2+k/\sigma}\quad\mbox{and}\quad r^*_{\varepsilon,k}\asymp\left(\varepsilon\left[\log {d\choose k}  \right]^{1/4}  \right)^{{4\sigma}/{(4\sigma+k)}}.
\end{gather*}
Therefore, under the assumption that
$\log {d\choose k}=\smallO \left( \log(\varepsilon^{-1}) \right)$, we obtain as $\varepsilon\to 0$
\begin{equation*}
	(\nu_{\lb}^*)^2 \asymp \varepsilon^{-2} r_{\varepsilon,k}^{2+k/\sigma}  \leq \varepsilon^{-2} (r_{\varepsilon,k}^*)^{2+k/\sigma} \asymp \varepsilon^{{2k}/{(4 \sigma + k)}} \left[ \log {d \choose k} \right]^{(2 \sigma +k)/{(4 \sigma + k)}} = \smallO(1).
\end{equation*}

Returning to (\ref{LB:risk_uneq1}), we may continue
\begin{align}
	R_{\varepsilon,k} &\geq \sum_{u \in \mathcal{U}_{k,d}} \inf_{\tilde{\etab}} \E{\pi_{\eta_u}} \operatorname{E}_{\pi,\eta_u} |\eta_u - \tilde{\eta}_u| = \sum_{u \in \mathcal{U}_{k,d}} \inf_{\tilde{\etab}} \left( (1-p) \E{\pi,0}(\tilde{\eta}_u) + p \E{\pi,1}(1-\tilde{\eta}_u) \right),
	\label{LB:risk_uneq2}
\end{align}
where $\tilde{\eta}_u$ may be viewed as a (nonrandomized) test and the quantity
\begin{equation*}
	\inf_{\tilde{\etab}} \left( (1-p) \E{\pi,0}(\tilde{\eta}_u) + p \E{\pi,1}(1-\tilde{\eta}_u) \right)
\end{equation*}
coincides with the Bayes risk in the problem of testing $H_0: \operatorname{P} = \operatorname{P}_{\pi,0}$ vs. $H_1: \operatorname{P} = \operatorname{P}_{\pi,1}$. The infimum over $\tilde{\eta}_u$ is attained for the Bayes test $\eta_B$ defined by (see, for example, Section 8.11 of \cite{DeGroot})
\begin{equation}
	\eta_B (\Yb_{\!\!u}) = \ind{\Lambda_{\pi,u} \geq \frac{1-p}{p}},
	\label{LB:opt_Bayes_test}
\end{equation}
where $\Lambda_{\pi,u}$ is the likelihood ratio defined in (\ref{LB:likelihood_ratio}).
It now follows from (\ref{LB:risk_uneq2}) and (\ref{LB:opt_Bayes_test}) that for any $u\in{\cal U}_{k,d}$ (from now on, we choose some $u$ and fix it)
\begin{align}
	R_{\varepsilon,k} &\geq {d \choose k} \left[ (1-p) \E{\pi,0} \left( 	\eta_B (\Yb_{\!\!u}) \right) + p \E{\pi,1} \left( 	\eta_B (\Yb_{\!\!u}) \right) \right] \nonumber \\
	&= {d \choose k}  (1-p) \operatorname{P}_{\pi,0} \left( \Lambda_{\pi,u} \geq  \frac{1-p}{p} \right) + {d \choose k} p \operatorname{P}_{\pi,1} \left( \Lambda_{\pi,u} < \frac{1-p}{p} \right) \nonumber \\
	& =: I_{\varepsilon,k}^{(1)} + I_{\varepsilon,k}^{(2)},
	\label{LB:risk_uneq3}
\end{align}
where both terms $I_{\varepsilon,k}^{(1)}$ and $I_{\varepsilon,k}^{(2)}$ are nonnegative.
Hence, it remains to show that at least one of these terms is positive for all small enough $\varepsilon$. To this end, we shall consider two (overlapping) cases
that cover assumptions (\ref{cond:inf2}) and (\ref{adas}).

\medskip
\textbf{Case 1:} Let $r_{\varepsilon,k}>0$ be such that
\begin{equation}
	0<\liminf_{\varepsilon \to 0} \frac{a_{\varepsilon,k}(r_{\varepsilon,k})}{\sqrt{\log {d\choose k}}} \leq\limsup_{\varepsilon \to 0} \frac{a_{\varepsilon,k}(r_{\varepsilon,k})}{\sqrt{\log {d\choose k}}} < \sqrt{2\beta}.
	\label{cond:inf2_lowpart}
\end{equation}
In this case, $ I_{\varepsilon,k}^{(2)}>0$  for all small enough $\varepsilon$.	
To show this, we introduce the random variables
\begin{align*}
	\lambda_{\pi} = \lambda_\pi (\Yb_{\!\!u})=\log \Lambda_{\pi,u} = \sum_{\lb \in \mathring{\mathbb{Z}}_u} \lambda_{\pi,\lb} (Y_\lb),
	\label{LB:def_lambda_pi}
\end{align*}
where
\begin{equation*}
	\lambda_{\pi,\lb}= \lambda_{\pi,\lb} (Y_\lb) =- \frac{(\nu_\lb^*)^2}{2} + \log \cosh (\nu_\lb^* Y_\lb),\quad \lb \in \mathring{\mathbb{Z}}_u,
	\label{LB:def_lambda_pi_lb}
\end{equation*}
with $\nu_\lb^*=\theta_{\lb}^*/\varepsilon$ obeying (\ref{LB:nu_star_is_o1}), and denote
\begin{equation}
	H = \log \frac{1-p}{p} = \log \left[ {d \choose k}^{\beta}-1\right].
	\label{LB:def_H}
\end{equation}
(We write $\lambda_{\pi}$ instead of $\lambda_{\pi,u}$ because $u$ has been chosen and fixed.)
Then, we can write
\begin{align}
	\operatorname{P}_{\pi,1} \left( \Lambda_{\pi,u} < \frac{1-p}{p} \right) &= \operatorname{P}_{\pi,1} \left( \lambda_\pi < H \right) = \E{\pi,1} \ind{\lambda_\pi < H} \nonumber \\
	&= \E{\pi,0} \left( e^{\lambda_\pi} \ind{\lambda_\pi < H} \right) = \E{\pi,0} \left( e^{h \lambda_\pi} e^{ \lambda_\pi (1-h)} \ind{\lambda_\pi < H} \right).
	\label{LB:P_pi_1_eq1}
\end{align}
Next, consider
the probability measures $\Po{h,\lb}, \; \lb \in \mathring{\mathbb{Z}}_u,$ defined by
\begin{equation}
	\frac{d \Po{h}}{d \Po{\pi,0}} (\Yb_{\!\!u}) = \frac{e^{h \sum_{\lb \in \mathring{\mathbb{Z}}_u} \lambda_{\pi,\lb}}}{\E{\pi,0} e^{h \sum_{\lb \in \mathring{\mathbb{Z}}_u} \lambda_{\pi,\lb} } } = \prod_{\lb \in \mathring{\mathbb{Z}}_u} \frac{e^{h \lambda_{\pi,\lb}}}{\E{\pi,0} e^{h \lambda_{\pi,\lb}} } =: \prod_{\lb \in \mathring{\mathbb{Z}}_u} \frac{d\Po{h,\lb}}{d \Po{\pi,0}} (Y_\lb).
	\label{LB:lem1_def_Phl}
\end{equation}
The following result holds true.

\begin{lemma} \label{lem:LB_lemma1}
	Let $h>0$ satisfy
	\begin{equation}
		\E{\Po{h}} (\lambda_\pi) = H.
		\label{LB:Elambda_is_H}
	\end{equation}
	Then, for $\lb \in \mathring{\mathbb{Z}}_u$  as $\varepsilon \to 0$
	\begin{align*}
		\E{\Po{h,\lb}} (\lambda_{\pi,\lb}) &= \frac{(\nu^*_\lb)^4}{2} \left(h-\frac{1}{2}\right) + \smallO((\nu^*_\lb)^4), \\ 
		\operatorname{var}_{\Po{h,\lb}}(\lambda_{\pi,\lb}) &= \frac{(\nu^*_\lb)^4}{2}  + \smallO((\nu^*_\lb)^4), 
	\end{align*}
	and hence
	\begin{align}
		\E{\Po{h}} (\lambda_{\pi}) &= a^2_{\varepsilon,u}(r_{\varepsilon,k}) \left(h-\frac{1}{2}\right)(1+\smallO(1)), \label{lem:LB_Elambda1}\\
		\sigma_h^2 := \operatorname{var}_{\Po{h}}(\lambda_{\pi}) &= a^2_{\varepsilon,u}(r_{\varepsilon,k}) (1+\smallO(1)). \label{lem:LB_varlambda1}
	\end{align}
\end{lemma}

In addition, we have the following result.

\begin{lemma} \label{lem:LB_lemma2}
	The Lyapunov condition for $(\lambda_{\pi,\lb})_{\lb \in \mathring{\mathbb{Z}}_u}$, $u\in {\cal U}_{k,d}$, is satisfied, that is,
	\begin{equation*}
		L_h^{(4)} :=\frac{\sum_{\lb \in \mathring{\mathbb{Z}}_u} \E{\Po{h,\lb}} (\lambda_{\pi,\lb} - \E{\Po{h,\lb}}(\lambda_{\pi,\lb}))^4}{\left( \sum_{\lb \in \mathring{\mathbb{Z}}_u} \var{\Po{h,\lb}} (\lambda_{\pi,\lb}) \right)^2} \to 0,	\quad \varepsilon\to 0,
	\end{equation*}
	and hence,  with respect to $\Po{h}$, the random variable $({\lambda_\pi - H})/{\sigma_h} $ has a standard normal distribution in the limit.
\end{lemma}

The proofs of these lemmas are similar to those of Lemmas 2 and 3 in \cite{ISt-2014}, and therefore are omitted.

We shall apply Lemma \ref{lem:LB_lemma1} to examine the right-hand side of (\ref{LB:P_pi_1_eq1}).
For this, we choose $h$ to satisfy (\ref{LB:Elambda_is_H}) and infer from (\ref{LB:def_H}) and (\ref{lem:LB_Elambda1}) that as $\varepsilon\to 0$
\begin{equation}
	\left( h-\frac{1}{2} \right) a_{\varepsilon,u}^2 (r_{\varepsilon,k}) (1+\smallO(1)) = H \sim \beta \log {d \choose k},
	\label{LB:H_is_ha2}
\end{equation}
which leads to
\begin{equation}
	h \sim \frac{1}{2} + \frac{\beta  \log {d \choose k}}{a_{\varepsilon,u}^2 (r_{\varepsilon,k})}.
	\label{LB:h_sim}
\end{equation}
In view of (\ref{cond:inf2_lowpart}), we have $a_{\varepsilon,u}^2 (r_{\varepsilon,k}) \asymp \log {d \choose k}$ as $\varepsilon\to 0$. From this and (\ref{LB:h_sim}),
\begin{equation}
	h = \mathcal{O}(1), \quad \varepsilon \to 0.
	\label{LB:h_bigO1}
\end{equation}

Omitting again the dependence on $u$, denote $\Psi (h) := \E{\pi,0} e^{h \lambda_\pi}$ and observe that
$$\Psi(h)= \prod_{\lb \in \mathring{\mathbb{Z}}_u} \E{\pi,0} e^{h \lambda_{\pi,\lb}} =: \prod_{\lb \in \mathring{\mathbb{Z}}_u} \Psi_\lb(h).$$ Then, in view of (\ref{LB:lem1_def_Phl}),
the probability in (\ref{LB:P_pi_1_eq1}) equals
\begin{align}
	\operatorname{P}_{\pi,1} \left( \Lambda_{\pi,u} < \frac{1-p}{p} \right) &= \E{\Po{h}} \left( \frac{d \Po{\pi,0}}{d\Po{h}} (\Yb_{\!\!u}) e^{h \lambda_\pi} e^{(1-h) \lambda_\pi} \ind{\lambda_\pi < H} \right) \nonumber \\
	&= \Psi (h) \E{\Po{h}} \left( e^{(1-h) \lambda_\pi} \ind{\lambda_\pi < H} \right).
	\label{LB:P_pi_1_eq2}
\end{align}
Define the functions
\begin{align*}
	p(v,y) = -\frac{v^2}{2} + \log \cosh (vy), \quad  \tilde{p} (v,y) = p(v,y)+ \frac{v^4}{4},
\end{align*}
and note that
\begin{equation*}
	\lambda_{\pi,\lb} = p(\nu^*_\lb,Y_\lb), \quad \lb \in \mathring{\mathbb{Z}}_u.
	\label{LB:lambda_pi_lb_is_p_fun}
\end{equation*}
It follows from relation (6.6) of \cite{BI-2011}
that for a positive value $h=\mathcal{O}(1)$ (see (\ref{LB:h_bigO1})) as $v\to 0$
\begin{equation*}
	\E{\Po{\pi,0}} e^{h \tilde{p}(v,Y_\lb)} = \exp \left( \frac{h^2 v^4}{4} + \smallO(v^4) \right).
\end{equation*}
Therefore, as $\varepsilon\to 0$
\begin{align*}
	\Psi_\lb(h) &= \E{\Po{\pi,0}} e^{h \lambda_{\pi,\lb}} = \E{\Po{\pi,0}} e^{h p(\nu_\lb^*,Y_\lb)} = \E{\Po{\pi,0}} e^{h \tilde{p}(\nu_\lb^*,Y_\lb)-h(\nu_\lb^*)^4/4} = \exp \left( \frac{h^2 - h}{4}(\nu_\lb^*)^4 + \smallO ((\nu_\lb^*)^4) \right),
\end{align*}
and hence, by using (\ref{def:a2_with_theta_star}) and the definition of $\nu_{\lb}^*$, we obtain
\begin{equation}
	\Psi(h) = \prod_{\lb \in \mathring{\mathbb{Z}}_u} \Psi_\lb(h) = \exp \left( \frac{h^2 - h}{2} a^2_{\varepsilon,u}(r_{\varepsilon,k}) (1+\smallO(1)) \right).
	\label{LB:Psi_by_h}
\end{equation}
Substituting (\ref{LB:Psi_by_h}) into (\ref{LB:P_pi_1_eq2}) gives
the following expression for the second term of the right-hand side of~(\ref{LB:risk_uneq3}):
\begin{align}
	I_{\varepsilon,k}^{(2)} &= {d \choose k} p \operatorname{P}_{\pi,1} \left( \Lambda_{\pi,u} < \frac{1-p}{p} \right)\nonumber\\ 
	& = {d \choose k}^{1-\beta}\exp\left( \frac{h^2-h}{2} a^2_{\varepsilon,u}(r_{\varepsilon,k}) (1+\smallO(1))  \right) \E{\Po{h}} \left( e^{(1-h) \lambda_\pi} \ind{\lambda_\pi < H} \right).
	\label{LB_I2_eq1}
\end{align}

Observe now that, under assumption (\ref{cond:inf2_lowpart}), by the choice of $h$ (see (\ref{LB:Elambda_is_H}), (\ref{lem:LB_Elambda1}) and (\ref{LB:h_sim})), as $\varepsilon\to 0$
\begin{align*}
	1-h \sim 1- \left( \frac{1}{2} + \frac{H}{a^2_{\varepsilon,u}(r_{\varepsilon,k})} \right) \sim \frac{1}{2} - \frac{\beta \log {d \choose k}}{a^2_{\varepsilon,u}(r_{\varepsilon,k})} < 0,
\end{align*}
that is, for all small enough $\varepsilon$,
\begin{equation*}
	1-h <0.
\end{equation*}
Having this property of $h$,
we can show that for a small $\varepsilon^{*}>0$,
the family of variables $\{\zeta_\varepsilon\}_{\varepsilon\in(0,\varepsilon^{*})}$, where $\zeta_\varepsilon:=e^{(1-h)\lambda_{\pi}}\mathds{1} (\lambda_\pi < H)$,
satisfies
\begin{gather}\label{ui}
	\sup_{0<\varepsilon<\varepsilon^{*}}\E{\Po{h}} \left\{ |\zeta_\varepsilon|\mathds{1} (|\zeta_\varepsilon|>N)\right\}\to 0, \quad N\to \infty,
\end{gather}
that is, $\{\zeta_\varepsilon\}_{\varepsilon\in(0,\varepsilon^{*})}$ is uniformly integrable.

Indeed, noting that $\zeta_\varepsilon\geq0$, we can write
\begin{gather}\label{A}
	\sup_{0<\varepsilon<\varepsilon^*}\E{\Po{h}} \left\{ |\zeta_\varepsilon|\mathds{1} (|\zeta_\varepsilon|>N)\right\}=
	\sup_{0<\varepsilon<\varepsilon^*}\E{\Po{h}} \left\{ e^{(1-h)\lambda_{\pi}}\mathds{1} \left(\lambda_{\pi}<-\frac{\log N}{h-1}\right)\right\},
\end{gather}
where, in view of Lemma \ref{lem:LB_lemma2},
\begin{gather}\label{B}
	\sup_{0<\varepsilon<\varepsilon^*}\Po{h}\left(\lambda_{\pi}<-\frac{\log N}{h-1}\right)\to 0,\quad N\to \infty.
\end{gather}
Next, by the definition of
$\lambda_{\pi}$, $\E{\operatorname{P}_{\pi,0}}(e^{\lambda_{\pi}})=1$, and hence
$\E{\Po{h}} \left( e^{(1-h)\lambda_{\pi}}\right)={\E{\operatorname{P}_{\pi,0}}(e^{\lambda_{\pi}})}/{{\Psi(h)}}=
{1}/{\Psi(h)},$
where, by means of~(\ref{LB:Psi_by_h}), $\Psi(h)\to \infty$ as $\varepsilon\to 0.$ Therefore,
for all $\varepsilon\in(0,\varepsilon^*)$, one can find a constant $C>0$ such that $1/\Psi(h)\leq C.$ This gives
\begin{gather}\label{C}
	\sup_{0<\varepsilon<\varepsilon^*}\E{\Po{h}} \left( e^{(1-h)\lambda_{\pi}}\right)<\infty.
\end{gather}
Combining (\ref{A}) to (\ref{C}) and using the absolute continuity of the Lebesgue integral, we arrive at (\ref{ui}).

Now, by making passage to the limit under the expectation sign and using Lemma \ref{lem:LB_lemma2}, we get as $\varepsilon\to 0$
\begin{align}
	\E{\Po{h}} \left[ e^{(1-h)\lambda_\pi} \ind{\lambda_\pi <H} \right] &\sim \int_{-\infty}^{H} \exp \left( (1-h)x \right) \frac{1}{\sqrt{2\pi}a_{\varepsilon,u}(r_{\varepsilon,k})} \exp\left( -\frac{(x-H)^2}{2 a^2_{\varepsilon,u}(r_{\varepsilon,k})} \right) dx \nonumber \\
	&=\frac{1}{\sqrt{2\pi}a_{\varepsilon,u}(r_{\varepsilon,k}) }  \exp\left((1-h)H+\frac{a^2_{\varepsilon,u}(r_{\varepsilon,k})(1-h)^2 }{2}  \right)\nonumber\\
	& \quad \times
	\int_{-\infty}^H \exp\left(-\frac{1}{2a^2_{\varepsilon,u}(r_{\varepsilon,k})}\left[x-(H+a^2_{\varepsilon,u}(r_{\varepsilon,k})(1-h))  \right]^2   \right)\, dx\nonumber\\
	&\geq  \exp\left((1-h)H+\frac{a^2_{\varepsilon,u}(r_{\varepsilon,k})(1-h)^2 }{2}  \right)
	\frac{1}{\sqrt{2\pi}a_{\varepsilon,u}(r_{\varepsilon,k}) }\nonumber\\ 
	&\quad \times \int_{-\infty}^{H+(1-h)a^2_{\varepsilon,u}(r_{\varepsilon,k})}
	\exp\left(-\frac{1}{2a^2_{\varepsilon,u}(r_{\varepsilon,k})}\left[x-(H+a^2_{\varepsilon,u}(r_{\varepsilon,k})(1-h))  \right]^2   \right)\, dx\nonumber\\
	&=\frac{1}{2}\exp\left((1-h)H+\frac{a^2_{\varepsilon,u}(r_{\varepsilon,k})(1-h)^2 }{2}  \right).
	\label{LB:Eeind_1_eq2}
\end{align}
From  (\ref{LB_I2_eq1}) and (\ref{LB:Eeind_1_eq2}), for all small enough $\varepsilon$
\begin{equation}
	I_{\varepsilon,k}^{(2)} \geq \frac{1}{2} {d \choose k}^{1-\beta}  \exp \bigg( \frac{h^2 - h}{2} a^2_{\varepsilon,u}(r_{\varepsilon,k}) (1+\smallO(1)) + H(1-h) +   \frac{a_{\varepsilon,u}^2(r_{\varepsilon,k}) (1-h)^2}{2} \bigg),
	\label{LB_I2_eq2}
\end{equation}
where, in view of (\ref{LB:H_is_ha2}),
\begin{gather*}
	\frac{h^2 - h}{2} a^2_{\varepsilon,u}(r_{\varepsilon,k}) (1+\smallO(1)) + H(1-h) + \frac{a_{\varepsilon,u}^2(r_{\varepsilon,k}) (1-h)^2}{2}
	= \smallO \left( a^2_{\varepsilon,u}(r_{\varepsilon,k}) \right) = \smallO \left( \log {d\choose k} \right).
\end{gather*}
From this and (\ref{LB_I2_eq2}), noting that $1-\beta>0$, we get
\begin{align*}
	I_{\varepsilon,k}^{(2)} &\geq \frac{1}{2} {d \choose k}^{1-\beta}  \exp \left( \smallO \left( \log {d\choose k} \right) \right) = \frac{1}{2} \exp \left( (1-\beta) \log {d \choose k} (1+\smallO(1)) \right)\nonumber \\
	& \geq \frac{1}{4} {d \choose k} ^{1-\beta} \to \infty, \quad \text{as }\varepsilon\to 0.
	%\label{LB_I2_eq3}
\end{align*}
Recalling (\ref{LB:risk_uneq3}), we now conclude that, under assumption (\ref{cond:inf2_lowpart}), for all small enough $\varepsilon$, the minimax risk $R_{\varepsilon,k}$ satisfies
\begin{equation*}
	R_{\varepsilon,k} \geq 	I_{\varepsilon,k}^{(2)} \geq {\rm const} >0.
\end{equation*}

\medskip
\textbf{Case 2:}  Let $r_{\varepsilon,k}>0$ be such that
\begin{equation}
	\sqrt{2}(1-\sqrt{1-\beta})<\liminf_{\varepsilon \to 0} \frac{a_{\varepsilon,k}(r_{\varepsilon,k})}{\sqrt{\log {d\choose k}}} \leq \limsup_{\varepsilon \to 0} \frac{a_{\varepsilon,k}(r_{\varepsilon,k})}{\sqrt{\log {d\choose k}}} < \sqrt{2} (1+\sqrt{1-\beta}).
	\label{cond:inf2_uppart}
\end{equation}
In this case, in view of (\ref{LB:risk_uneq3}), it is sufficient to show  that for some chosen $u\in {\cal U}_{k,d}$ and all small enough~$\varepsilon$
\begin{equation}
	I_{\varepsilon,k}^{(1)} = {d \choose k} (1-p) \Po{\pi,0} \left( \Lambda_{\pi,u} \geq \frac{1-p}{p} \right) >0,
	\label{LB:I1_positive}
\end{equation}
The proof of inequality (\ref{LB:I1_positive}), as presented below, is similar to that of (34) in \cite{ISt-2014}.

First, using the same notation as in Case 1, we can write
\begin{align*}
	& \operatorname{P}_{\pi,0} \left( \Lambda_{\pi,u} \geq \frac{1-p}{p} \right) = \operatorname{P}_{\pi,0} \left( \lambda_\pi \geq H \right) \nonumber \\ 
	&\quad = \E{\pi,0} \ind{\lambda_\pi \geq H}  
	= \E{\pi,0} \left( e^{h \lambda_\pi} e^{ -h\lambda_\pi} \ind{\lambda_\pi \geq H} \right)\nonumber\\
	&\quad =\E{\Po{h}} \left( \frac{d \Po{\pi,0}}{d\Po{h}} (\Yb_{\!\!u}) e^{h \lambda_\pi} e^{-h \lambda_\pi} \ind{\lambda_\pi \geq H} \right) = \Psi (h) \E{\Po{h}} \left( e^{-h \lambda_\pi} \ind{\lambda_\pi \geq H} \right) \nonumber \\
	& \quad = \Psi (h) e^{-h H} \E{\Po{h}} \left( e^{-h (\lambda_\pi-H)} \ind{\lambda_\pi-H \geq 0} \right) =: \Psi (h) e^{-h H} J(h), 
\end{align*}
which gives
\begin{equation}
	\log \Po{\pi,0} \left( \Lambda_{\pi,u} \geq \frac{1-p}{p} \right) = \log \Psi(h) - hH + \log J(h).
	\label{LB:logP0}
\end{equation}
Next, choose again $h$ as in (\ref{LB:Elambda_is_H}), denote
$Z_h := ({\lambda_{\pi} - H})/{\sigma_h},$ where $\sigma_h^2 = \var{\Po{h}} (\lambda_{\pi}),$
and observe that $\E{\Po{h}}(Z_h) = 0$ and $\var{\Po{h}}(Z_h)=1$.
By the definition of $J(h)$, for some positive $\tau = \tau_\varepsilon \to 0$
\begin{align}
	J(h) & = \E{\Po{h}} \left( e^{-h (\lambda_\pi-H)} \ind{\lambda_\pi-H \geq 0} \right) = \E{\Po{h}} \left( e^{-h Z_h \sigma_h} \ind{Z_h \sigma_h \geq 0} \right) \nonumber \\
	&\geq \E{\Po{h}} \left( e^{-h Z_h \sigma_h} \ind{Z_h \in [0,\tau]} \right) \geq  e^{-h \tau \sigma_h} \Po{h} \left( Z_h \in [0,\tau] \right). \nonumber
\end{align}
From this, for all small enough $\varepsilon$, we obtain
\begin{equation}
	e^{-h \tau \sigma_h} \Po{h} \left( Z_h \in [0,\tau] \right) \leq J(h) \leq 1,
	\label{LB:boundary_for_Jh}
\end{equation}
where, by Lemma \ref{lem:LB_lemma2},
\begin{equation*}
	\Po{h} \left( Z_h \in [0,\tau] \right) \sim \int_{0}^{\tau} \frac{1}{\sqrt{2 \pi}} e^{-t^2/2} dt \sim \frac{\tau}{\sqrt{2 \pi}} >0,
\end{equation*}
and hence (\ref{LB:boundary_for_Jh}) gives
\begin{equation}
	-\tau h \sigma_h + \log \tau - \frac{1}{2} \log (2\pi) \leq \log J(h) \leq 0.
	\label{ineq1}
\end{equation}
Next, in view of relations (\ref{LB:def_H}), (\ref{lem:LB_varlambda1}), (\ref{LB:H_is_ha2}) and (\ref{LB:h_bigO1}), which are also valid in Case~2 under study, we have as $\varepsilon\to 0$
\begin{equation}
	h \sigma_h \asymp a_{\varepsilon,u} (r_{\varepsilon,k}) \asymp \sqrt{H} \asymp \sqrt{\log {d \choose k}} .
	\label{LB:hsigma_asymp_sqrtH}
\end{equation}
which gives
\begin{equation*}
	h \sigma_h \to \infty.
\end{equation*}
If, in addition, we choose $\tau = \tau_\varepsilon \to 0$ to have
$\tau^{-1}\log(\tau^{-1})=\smallO( h\sigma_h)$ as $\varepsilon\to 0,$
we get from (\ref{ineq1}) and (\ref{LB:hsigma_asymp_sqrtH}) that
\begin{equation*}
	\log J(h) = \smallO (h\sigma_h) = \smallO(\sqrt{H}), \quad \varepsilon\to 0,
\end{equation*}
and hence, recalling  (\ref{LB:h_bigO1}) and (\ref{LB:Psi_by_h}), we obtain from (\ref{LB:logP0})
\begin{eqnarray*}
	\log \Po{\pi,0} \left( \Lambda_{\pi,u} \geq \frac{1-p}{p} \right) &=& \log \Psi(h) - hH + \smallO (\sqrt{H})
	\nonumber\\	&=& \frac{h^2-h}{2}  a_{\varepsilon,u}^2 (r_{\varepsilon,k}) +\smallO(H) - hH + \smallO(\sqrt{H}).
\end{eqnarray*}
Now, expressing $h$ as (see (\ref{LB:H_is_ha2}))
$	h= {1}/{2} + {H}/{a_{\varepsilon,u}^2 (r_{\varepsilon,k})} (1+\smallO(1)),$
we may continue
\begin{gather}
	\log \Po{\pi,0} \left( \Lambda_{\pi,u} \geq \frac{1-p}{p} \right)
	= -\frac{1}{2 a_{\varepsilon,u}^2 (r_{\varepsilon,k})} \left( H + \frac{a_{\varepsilon,u}^2 (r_{\varepsilon,k})}{2} \right)^2 (1+\smallO(1)).
	\label{LB:logP0_eq4}
\end{gather}
Exponentiating both sides of  (\ref{LB:logP0_eq4}) gives as $\varepsilon\to 0$
\begin{align*}
	\Po{\pi,0} \left( \Lambda_{\pi,u} \geq \frac{1-p}{p} \right) &= \exp \left\{ -\frac{1}{2 a_{\varepsilon,u}^2 (r_{\varepsilon,k})} \left( H + \frac{a_{\varepsilon,u}^2 (r_{\varepsilon,k})}{2} \right)^2 (1+\smallO(1)) \right\} 
\end{align*}
and hence
\begin{equation*}
	I_{\varepsilon,k}^{(1)}=	
	\exp \left\{\left(\log {d\choose k} -\frac{1}{2 a_{\varepsilon,u}^2 (r_{\varepsilon,k})} \left( H + \frac{a_{\varepsilon,u}^2 (r_{\varepsilon,k})}{2}\right)^2  \right) (1+\smallO(1))\right\}.
\end{equation*}
It remains to show that for all sufficiently small $\varepsilon$
\begin{equation*}
	\log {d\choose k} - \frac{\left( H + {a_{\varepsilon,u}^2 (r_{\varepsilon,k})}/{2} \right)^2}{2 a_{\varepsilon,u}^2 (r_{\varepsilon,k})} >0.
\end{equation*}
The latter inequality holds true when, for all small enough $\varepsilon$, it holds that
$	\sqrt{2 \log {d\choose k}} \left( 1 - \sqrt{1-\beta} \right) < a_{\varepsilon,u}^2 (r_{\varepsilon,k}) < \sqrt{2 \log {d\choose k}} \left( 1 + \sqrt{1-\beta} \right)$,
which, due to~(\ref{LB:def_H}), is guaranteed by (\ref{cond:inf2_uppart}). The validity of inequality (\ref{LB:I1_positive}) is thus verified and, in view of
(\ref{LB:risk_uneq3}), for all small enough~$\varepsilon$
\begin{equation*}
	R_{\varepsilon,k} \geq 	I_{\varepsilon,k}^{(1)} \geq {\rm const} >0.
\end{equation*}

Combining Cases 1 and 2 and noting that $1-\sqrt{1-\beta}<\sqrt{\beta}$ for all $\beta\in(0,1)$ completes the proof of Theorem \ref{theorem:LB_figed_k}. \qedsymbol
\medskip

\textit{Proof of Theorem \ref{th:UB_growingK}.} When $d=d_\varepsilon\to \infty$,  $k=k_\varepsilon\to \infty$, and $k=\smallO(d)$ as $\varepsilon\to 0$,
we have
\begin{gather}\label{ausl}
	{d \choose k}\sim {d^k}/{k!}\quad\mbox{and}\quad \log {d \choose k}\sim k\log\left({d}/{k}\right).
\end{gather}
In particular, it is still true that ${d \choose k}\to \infty$ as $\varepsilon\to 0$. Therefore, the proof largely repeats that of Theorem \ref{th:UB_fixedK},
and we only highlight the differences.

Due to (\ref{ausl}), the assumption $\log {d \choose k}=\smallO(\log \varepsilon^{-1})$ ensures that $k=\smallO(\log \varepsilon^{-1})$,
which, in view of the assumption
$\log\log d=\smallO(k)$, is the same as
\begin{gather}\label{usl}
	k=\smallO\left(\log \left(\varepsilon\left[\log {d\choose k}\right]^{1/4} \right)^{-1}  \right).
\end{gather}
With the restriction on $k$ as in (\ref{usl}),  by means of
(\ref{def:r_star}), (\ref{aek1}), and (\ref{def:r_star_ekm}), as $\varepsilon\to 0$, cf. (\ref{rek}),
\begin{gather}\label{rek1}
	r_{\varepsilon,k}^*\asymp r_{\varepsilon,k,m}^*\asymp\left(\varepsilon\left[\log {d\choose k}\right]^{1/4} \right)^{{4\sigma}/{(4\sigma+k)}} k^{-{\sigma}/{2}},
	\quad m=1,\ldots,M_k.
\end{gather}
(For more details on the derivation of (\ref{rek1}) and the importance of condition (\ref{usl}), we refer to Section~4 of \cite{IS-2015}.)
In particular,  $r_{\varepsilon,k}^*=\smallO(k^{-{\sigma}/{2}})$ as $\varepsilon\to 0$, as required:
the hypotheses $\mathbf{H}_{0,u}: \thetab_u=\boldsymbol{0} $ and $\mathbf{H}^{\varepsilon}_{1,u}: \thetab_u \in \mathring{\Theta}_{\cb_u}(r_{\varepsilon,k})$ separate asymptotically when
$r_{\varepsilon,k}/r^*_{\varepsilon,k}\to \infty$, 
and the range of interesting values of $r_{\varepsilon,k}$ is $r_{\varepsilon,k}\in(0,(2\pi)^{-\sigma}k^{-{\sigma}/{2}})$, as explained in Section \ref{PS}.

Next, for all indices $\lb \in \mathring{\mathbb{Z}}^k$ with the property $\theta^*_\lb(r_{\varepsilon,k})\neq 0$,  
in view of (\ref{def:omega}), (\ref{theta2}), and (\ref{aek1}), as $\varepsilon\to 0$
\begin{gather}\label{omega8}
	\omega_\lb(r_{\varepsilon,k})=
	\frac{[\theta^*_\lb(r_{\varepsilon,k})]^2}{2\varepsilon^2 a_{\varepsilon,k}(r_{\varepsilon,k})}
	\asymp r_{\varepsilon,k}^{{k}/{(2\sigma)}}\left({2\pi k}/{e}\right)^{{k}/{4}} k^{5/4}.
\end{gather}
Combining relations (\ref{rek1}) and (\ref{omega8}), we obtain for all $m=1,\ldots,M_k$, as $\varepsilon\to 0$, cf. (\ref{omegamax}),
\begin{gather}\label{omega9}
	\max_{\lb \in \mathring{\mathbb{Z}}^k} \omega_\lb(r^*_{\varepsilon,k,m})\asymp
	\left( \varepsilon \left[ \log {d \choose k} \right]^{1/4} \right)^{{2k}/{(4 \sigma + k)}}\left({2\pi}/{e}\right)^{{k}/{4}} k^{5/4}.
\end{gather}

The asymptotic expression in (\ref{omega9}) differs from that in (\ref{omegamax}) by a factor tending to infinity.
Nevertheless, thanks to the assumptions on $d$ and $k$,
the rest of the proof goes along the lines of that of Theorem \ref{th:UB_fixedK}
and therefore is omitted. \qedsymbol

\medskip
\textit{Proof of Theorem \ref{theorem:LW_growing_k}.}
Under the assumed conditions on $d$ and $k$, it is true that ${d \choose k}\to \infty$ as $\varepsilon\to 0$.
Hence, the proof is similar to that of Theorem \ref{theorem:LB_figed_k}.
It is only required to show that, for some chosen $u\in {\cal U}_{k,d}$, when $r_{\varepsilon,k}>0$ is such that  $\limsup_{\varepsilon\to 0}r_{\varepsilon,k}/r^*_{\varepsilon,k}<1$ with
$r^*_{\varepsilon,k}$ as in (\ref{rek1}), the quantities $\nu_{\lb}^*=\theta^*_{\lb}(r_{\varepsilon,k})/\varepsilon$, $\lb\in  \mathring{\mathbb{Z}}_u$, satisfy,
cf. (\ref{LB:nu_star_is_o1}),
\begin{gather*}
	\nu_{\lb}^*=\smallO(1),\quad \varepsilon\to 0.
\end{gather*}

For this, we first note that as $\varepsilon\to 0$
\begin{gather*}
	\log {d \choose k}=\smallO(\log \varepsilon^{-1})\quad\mbox{and}\quad k=\smallO(\log \varepsilon^{-1}),
\end{gather*}
Then, using (\ref{theta2}) and (\ref{rek1}), for all indices $\lb \in \mathring{\mathbb{Z}}_u$ with the property $\theta^*_\lb(r_{\varepsilon,k})\neq 0$, we have, as $\varepsilon\to 0$,
\begin{eqnarray*}
	\left[\nu_{\lb}^*\right]^2&=&\varepsilon^{-2}\left[\theta^*_{\lb}(r_{\varepsilon,k})\right]^2
	\asymp\varepsilon^{-2}r_{\varepsilon,k}^{2+k/\sigma}\left({2\pi k}/{e}\right)^{{k}/2} k^{3/2} \leq \varepsilon^{-2}(r^*_{\varepsilon,k})^{2+k/\sigma}\left({2\pi k}/{e}\right)^{{k}/2} k^{3/2}\\
	&=& \varepsilon^{2k/(4\sigma+k)}\left[\log {d\choose k}  \right]^{(2\sigma+k)/(4\sigma+k)}
	\left({2\pi}/{e}\right)^{{k}/2} k^{-\sigma+3/2}=\smallO(1).
\end{eqnarray*}

\qedsymbol

\section*{Acknowledgements} 
The research of N. Stepanova was supported by an NSERC grant. The research of M. Turcicova was partially supported by NSERC grants during the author's stays at Carleton University in 2022 and 2023.

%% or include bibliography directly:


\begin{thebibliography}{9}

	\bibitem{BI-2011}
\textsc{Butucea, C.} and \textsc{Ingster, Yu.} (2011): {Detection of a sparse submatrix of a high-dimensional noisy
	matrix}, \textit{https://arxiv.org/abs/1109.0898v1}.

\bibitem{BI-2012}
\textsc{Butucea, C.} and \textsc{Ingster, Yu. (2013)}: {Detection of a sparse submatrix of a high-dimensional noisy
	matrix}, \textit{Bernoulli} \textbf{19}(5B), 2652--2688.
%	\MR{3160567}

\bibitem{BSt-2017}
\textsc{Butucea, C.} and \textsc{Stepanova, N.} (2017): Adaptive variable selection in nonparametric sparse models, \textit{Electronic Journal of Statistics} \textbf{11}, 2321--2357.
%\MR{3656494}

\bibitem{DeGroot}
\textsc{DeGroot, M.} (1970). \textit{Optimal Statistical Decisions}, McGraw-Hill Book Company, New York.
%\MR{0356303}

\bibitem{GS-2020} 
\textsc{Gao, Z.} and \textsc{Stoev, S.} (2020): {Fundamental limits of exact support recovery in high dimensions}, \textit{Bernoulli} \textbf{26}(4), 2605--2638.
%\MR{4140523}

\bibitem{DSA-2012}
\textsc{Gayraud, G.} and \textsc{Ingster, Yu.} (2012): {Detection of sparse additive functions}, \textit{Electronic Journal of Statistics} \textbf{6}, 1409--1448. 
%\MR{2988453}

\bibitem{Genovese-2012}
\textsc{Genovese, C.R., Jin, J., Wasserman, L.} and \textsc{Yao, Z.} (2012): {A comparison of the lasso and marginal regression}, \textit{Journal of Machine Learning Research} \textbf{13}, 2107--2143. 
%\MR{2956354}

\bibitem{GN2016} 
\textsc{Gin\'{e}, E.} and \textsc{Nickl, R.} (2016).
\textit{Mathematical Foundations of Infinite-Dimensional Statistical Models}, Cambridge University Press, New York.
%\MR{3588285}

\bibitem{IH.97} 
\textsc{Ibragimov, I. A.} and \textsc{Khasminskii, R. Z.} (1997). Some estimation
problems on infinite dimensional Gaussian white noise. In: {\it
Festschrift for Lucien Le Cam. Research papers in Probability and Statistics}, 275--296. Springer--Verlag, New York.
%\MR{1462950}

\bibitem{I1993a}
\textsc{Ingster, Yu. I.} (1993):
{Asymptotically minimax hypothesis testing for nonparametric alternatives. I}.
\textit{Mathematical Methods of Statistics} \textbf{2}(2), 85--114.
%\MR{1257978}

\bibitem{I1993b}
\textsc{Ingster, Yu. I.} (1993):
{Asymptotically minimax hypothesis testing for nonparametric alternatives. II}.
\textit{Mathematical Methods of Statistics} \textbf{2}(3), 171--189.
%\MR{1257983}

\bibitem{I1993c}
\textsc{Ingster, Yu. I.} (1993):
{Asymptotically minimax hypothesis testing for nonparametric alternatives. III}.
\textit{Mathematical Methods of Statistics} \textbf{2}(4), 249--268.
%\MR{1259685}

\bibitem{NGF-2003}
\textsc{Ingster, Yu.} and \textsc{Suslina, I.} (2003): \textit{Nonparametric Goodness-of-Fit Testing under Gaussian Models}, Springer, New York. 
%\MR{1991446}

\bibitem{IS-2005}
\textsc{Ingster, Yu.} and \textsc{Suslina, I.} (2005): {On estimation and detection of smooth
	functions of many variables}, \textit{Mathematical Metods of Statistics} \textbf{14}(3), 299--331.
%\MR{2195328}

\bibitem{ISt-2014}
\textsc{Ingster, Yu.} and \textsc{Stepanova, N.} (2014): {Adaptive variable selection in nonparametric sparse regression},  \textit{Journal of Mathematical Sciences} \textbf{199}(2), 184--201. 
%\MR{3032218}

\bibitem{IS-2015}
\textsc{Ingster, Yu.} and \textsc{Suslina, I.} (2015): {Detection of a sparse variable functions}, \textit{Journal of Mathematical Sciences} \textbf{206}(2), 181--196.
%\MR{3373876}

\bibitem{Ji-2012}
\textsc{Ji, O.} and \textsc{Jin, J.} (2012): {UPS delivers optimal phase diagram in high-dimensional variable selection}, \textit{Annals of Statistics} \textbf{40}(1),  73--103.
%\MR{3013180}

\bibitem{Lin-2000} 
\textsc{Lin, Y.} (2000). Tensor product space ANOVA models, \textit{Annals of Statistics} \textbf{28}, 734--755.
%\MR{1792785}

\bibitem{MSt-2023} 
\textsc{Miller, J. C.} and \textsc{Stepanova, N. A.} (2023): Adaptive signal recovery with Subbotin noise,
\textit{Statistics and Probability Letters} \textbf{196}(C), article 109791.
%\MR{4549624}

\bibitem{SK} 
\textsc{Skorohod, A. V.}  (1974). {\it Integration in Hilbert Spaces}. Springer--Verlag,
Berlin--New York.
%\MR{0466482}

\end{thebibliography}
\end{document}